\documentclass[10pt,a4paper,draft]{amsart}
\usepackage{amsmath}
\usepackage{amsthm}
\usepackage{amsfonts}
\usepackage{amssymb}
\usepackage{a4,a4wide}
\usepackage{leqno}
\usepackage[english]{babel}

\usepackage{graphics}
\usepackage{euscript}

\usepackage{color}
\usepackage[T1]{fontenc}
\usepackage{palatino}

\newtheorem{theorem}{Theorem}[section]
\newtheorem{lemma}[theorem]{Lemma}

\newtheorem{claim}[theorem]{Claim}

\theoremstyle{definition}

\theoremstyle{remark}
\newtheorem{remark}[theorem]{Remark}

\numberwithin{equation}{section}
\def\O{{\Omega}}
\def\o{{\omega}}

\def\eps{{\epsilon}}

\def\f{{\mathcal{F}}}
\def\h{{\mathcal{H}}}
\def\I{{\mathcal{I}}}
\def\j{{\mathcal{J}}}

\def\r{{\mathcal{R}}}

\def\D{{\mathcal{D}}}

\def\R{{\mathbb{R}}}

\def\N{{\mathbb{N}}}

\newcommand{\oph}[3]{\h_{_{{#1},{#2}}}[{#3}]}

\newcommand{\dem}[1]{\vskip 0.2\baselineskip \noindent {\bf{#1}}\vskip 0.2\baselineskip }
\newcommand{\fdem}{\vskip 0.2 pt \hfill $\square$ }

\newcommand{\nl}[2]{\|{#1}\|_{L^2{(#2)}}}
\newcommand{\nlto}[1]{\|{#1}\|_{2}}
\newcommand{\nlmu}[2]{\|{#1}\|_{L_\mu^2{(#2)}}}
\newcommand{\nlp}[3]{\|{#1}\|_{L^{#2}{(#3)}}}

\def\tilde{\widetilde}

\author{Jerome Coville}
\address{\noindent J. Coville -- INRA PACA, Equipe BIOSP, Centre de Recherche d'Avignon, Domaine Saint
Paul, Site Agroparc, 84914 Avignon cedex 9, France}
\email{jerome.coville@avignon.inra.fr}
\title{Convergence to equilibrium for positive solutions of some mutation-selection model }
\date\today
\begin{document}
\maketitle
\begin{abstract}
In this paper we are interested in  the long time behaviour of the positive solutions of the  mutation selection model with Neumann Boundary condition: 
$$
\frac{\partial u(x,t)}{dt}=u\left[r(x)-\int_{\O}K(x,y)|u|^{p}(y)\,dy\right]+\nabla\cdot\left(A(x)\nabla u(x)\right),\qquad \text{ in }\quad \R^+\times\O$$
where $\O\subset \R^N$ is a bounded smooth domain, $k(.,.) \in C(\bar \O \times C(\bar\O), \R), p\ge 1$ and $A(x)$ is a smooth elliptic matrix.
 
 In a blind competition situation, i.e $K(x,y)=k(y)$,  we show the existence of a unique positive  steady  state which is  positively globally stable.  That is, the positive steady state  attracts all the possible trajectories  initiated from any non negative initial datum. When $K$ is a general positive  kernel, we also present a necessary and sufficient condition  for  the existence of a positive steady states. We prove also some stability result on the dynamic  of the equation   when  the competition  kernel $K$ is of the form $K(x,y)=k_0(y)+\eps k_1(x,y)$. That is, we prove that for sufficiently small $\eps$  there exists a unique   steady state, which in addition is positively asymptotically stable.    
  The proofs of the global stability of the steady state essentially rely on  non-linear relative entropy  identities  and an orthogonal decomposition. These identities combined with the decomposition provide us some a priori estimates and differential inequalities essential to characterise  the asymptotic behaviour of the solutions. 
\end{abstract}

\section{Introduction and Main results}\label{s:intro}
In this paper we are interested in  the long time behaviour of the positive solutions of the nonlocal equation 
\begin{align}\label{msedp-eq-intro}
&\frac{\partial u(t,x)}{\partial t}=u(t,x)\left[r(x)-\int_{\O}K(x,y)|u(t,y)|^{p}\,dy\right]+ \nabla\cdot\left(A(x)\nabla u(t,x)\right)\qquad \text{ in }\quad \R^+\times\O\\
&\frac{\partial u(t,x)}{\partial n}=0, \quad \text{ in } \quad \R^+\times\partial\O\label{msedp-eq-intro-bc}\\
&u(0,x)=u_0(x) \label{msedp-eq-intro-ci}
\end{align}
where $\O\subset \R^N$ is a bounded smooth domain, $r(x)\in C^{0,1}(\bar \O)$ is positive, $p\ge 1$, $K(.,.) \in C^{0,1}(\bar \O \times \bar \O)$ and $A(x) \in \mathcal{M}_{n\times n}(\R)$ is  a uniform smooth ($C^{1,\alpha}$) elliptic matrix.

Such type of nonlocal model has been introduced to capture the evolution of a  population structured by a phenotypical trait \cite{Burger2000,Burger1994,Dieckmann2004,Perthame2007}.   In this context  $u(x,t)$ represents the density of a population at the phenotypical trait $x$  at time $t$,  which is submitted to two essential interactions:  mutation and selection.  
Here, the mutation process, which acts  as a diffusion operator on the traits space, is modelled by a classical diffusion operator whereas the selection process is modelled by the nonlocal term $u(t,x)\int_{\O}K(x,y)|u(t,y)|^{p}\,dy$.  In the literature, the selection  operator takes often the form $u(t,x)\int_{\O}K(x,y)|u(t,y)|\,dy$ \cite{Barles2008,Burger2000,Perthame2007}. 
A rigorous derivation of these equations  from stochastic processes can be found in \cite{Champagnat2008,Fournier2004}.

To our knowledge, a large part of the analysis of the long time behaviour of solutions of \eqref{msedp-eq-gen} concerns either situations where no mutation occurs  \cite{Barles2009,Barles2008,Burger2000,Calsina2005,Calsina2007,Carrillo2007,
Desvillettes2008,Diekmann2005,Jabin2011} or in the context of "adaptive dynamics", i.e. the evolution of the population is driven by small mutations,   \cite{Burger2000,Burger1994,Canizo,Carrillo2007,Champagnat2011, Lorz2011} and references therein.

In the latter case, the matrix $A(x):=\eps A_0(x)$ and   some asymptotic regimes are studied when $\eps \to 0$.   In this situation, an extensive work have been done  in developing a constrained Hamilton-Jacobi approach in order to analyse the long time behaviour of positive solutions of this type of models see for instance \cite{Barles2009,Barles2008,Carrillo2007,Champagnat2011, Diekmann2005}.  

Analysis of variants of \eqref{msedp-eq-intro} involving a nonlocal mutation process of the form \linebreak $\eps \int_{\O}\mu(x,y)(u(t,y)-u(t,x))\,dy$ instead of an elliptic diffusion can be found \cite{Calsina2005,Calsina2007,Calsina2012,Raoul2011,Raoul2012}. For these variants,  approaches based on semi-group theory  have been developed  to analyse the asymptotic behaviour and local stability of the positive stationary solution of \eqref{msedp-eq-intro} when $\eps \to 0$,   see  \cite{Calsina2005,Calsina2007,Calsina2012}.

In all those works, the small mutation assumptions appears to be a key feature in the analysis. 
 Our goal here is to analyse the long time behaviour of the solution to \eqref{msedp-eq-intro} -- \eqref{msedp-eq-intro-ci} in   situations where no restriction on the mutation operator are imposed. In particular, we want to understand situations where the rate of mutations is not  small compared to selection. This appears for example  in some virus population where the rate of mutation per reproduction cycle is high \cite{Cuevas2005,Fabre2012,Sanjuan2010,Zhu2009}.

In what follows, we will always make the following assumptions on $r$, $K$

\begin{equation}
\label{msedp-hyp1}
\left\{
\begin{aligned}
&\text{$A \in  \mathcal{M}_{n\times n}(\R)$ is a smooth uniform elliptic matrix,}\\
&\text{$r \in C^{0,1}(\O)$ is positive,}\\
&\text{$\O$ is a bounded Lipschitz domain in $\R^N$.  }\\
&\text{ $K \in C^{0,1}(\bar \O\times \bar \O), K> 0, $}
\end{aligned}
\right.
\end{equation}

Under the above assumptions the existence of a positive solution to the Cauchy problem \eqref{msedp-eq-intro}--\eqref{msedp-eq-intro-ci} is guarantee. Namely, we can easily prove
\begin{theorem}\label{msedp-thm-exits}
Assume $A,r,K$ satisfy \eqref{msedp-hyp1} and  $p\ge 1$ then for all $u_0 \in L^p(\O)$ there exists a positive smooth solution $u$  to  \eqref{msedp-eq-intro} -- \eqref{msedp-eq-intro-ci} so that  $u\in C([0,+\infty),L^p(\O))\cap C^1((0,+\infty),C^{2,\alpha}(\O))$.
\end{theorem}

The main problematic then remains to characterise the long time   behaviour of these solutions. In this direction our first result concerns the situations of blind competition, that is when the kernel $K(x,y)$ is independent of $x$. In this context the equations \eqref{msedp-eq-intro} -- \eqref{msedp-eq-intro-ci} rewrite
 \begin{align}
&\frac{\partial u}{\partial t}(t,x)=u(t,x)\left(r(x)-\int_{\O}k(y)|u(t,y)|^p\,dy\right)+\nabla\cdot(A(x)\nabla u(t,x) )\quad \text{ in } \quad \R^+\times\O \label{msedp-eq-red}\\
&\frac{\partial u}{\partial n}(t,x)=0 \quad \text{ in } \quad \R^+\times\partial\O\label{msedp-eq-red-bc}\\
&u(x,0)=u_0(x)\quad \text{ in } \quad \O.\label{msedp-eq-red-ci}
\end{align}

In this situation, we have 
\begin{theorem}\label{msedp-thm1}
Assume $A,r,k$ satisfy \eqref{msedp-hyp1} and  $p\ge 1$.
 Let $\lambda_1$ be the first eigenvalue  of the operator $ \nabla\cdot(A(x)\nabla  )+r(x)$ with Neumann boundary condition and let $\phi_1$ be a positive eigenfunction associated with $\lambda_1$, that is $\phi_1$ satisfies 
\begin{align}
&\nabla\cdot(A(x)\nabla \phi_1 )+r(x)\phi_1=-\lambda_1\phi_1 \quad \text{ in } \quad \O,\label{msedp-eq-pev}\\
&\frac{\partial \phi_1}{\partial n}(x)=0\quad \text{ on } \quad \partial\O.
\end{align}
Then we have the following asymptotic behaviour for any positive smooth ( at least $C^2$) solution $u(t,x)$ to \eqref{msedp-eq-red} -- \eqref{msedp-eq-red-bc}
\begin{itemize}
\item if $\lambda_1\ge 0$, there is no positive stationary solution and $u(t,x)\to 0$ as $t\to \infty$ 
\item if $\lambda_1< 0$, then $$u(t,x) \to \mu \phi_1 $$ where 
$\mu=\left(\frac{-\lambda_1}{\int_{\O}k(y)|\phi_1|^p(y)\,dy}\right)^{\frac{1}{p}}$ and $\phi_1$ has been normalized by $\nl{\phi_1}{\O}=1$. 
\end{itemize}
\end{theorem}

Next we establish an optimal  existence criteria for the positive stationary solution to \eqref{msedp-eq-intro}-\eqref{msedp-eq-intro-bc}. Namely, we prove
  
 \begin{theorem}\label{msedp-thm2}
Assume $A,r,K$ satisfy \eqref{msedp-hyp1} and  $p\ge 1$. Then there exists at least a positive smooth solution $\bar u$ of \eqref{msedp-eq-intro} -- \eqref{msedp-eq-intro-ci} \textbf{if and only if } $\lambda_1< 0$,  where $\lambda_1$ is defined  in Theorem \ref{msedp-thm1}. 
\end{theorem}

 Finally, we prove that the dynamic observed for blind selection kernel $K(x,y)=k(y)$ still holds for some perturbation of $k$. More precisely, let us consider a kernel $k_\eps(x,y)=k_0(y)+\eps k_1(x,y)$ with $k_i$ satisfying the assumption \eqref{msedp-hyp1}, then 
 we have the following
  \begin{theorem}\label{msedp-thm3}
Assume $A,r,K$ satisfy \eqref{msedp-hyp1} and  $p=1$ or $p=2$. Assume further that $K=k_\eps$ and let $u(t,x)$ be a positive smooth solution to \eqref{msedp-eq-intro}--\eqref{msedp-eq-intro-bc} with $K=k_\eps$. Then we have the following asymptotic behaviour: 
\begin{itemize}
\item if $\lambda_1\ge 0$, there is no positive stationary solution and  $u(t,x)\to 0$ as $t\to \infty$ uniformly. 
\item if $\lambda_1< 0$, then there exists $\eps^*$ so that for all $\eps \le \eps^*$ there exists a unique positive globally attractive equilibrium $\bar u_\eps$ to \eqref{msedp-eq-intro}-\eqref{msedp-eq-intro-bc} i.e. for all $u_0\ge_{\not \equiv} 0$, then   we have for all $x\in \O,$  
$$\lim_{t\to\infty}u(t,x) \to \bar u_\eps (x). $$  
\end{itemize}
\end{theorem}

 \subsection{Comments}
Before going to the proofs of these results, we would like to make some comments.
First, it comes directly from the proofs that  the Theorems \ref{msedp-thm1} and \ref{msedp-thm2} can be generalised to more general selection process. In particular, Theorem \ref{msedp-thm1} holds true if instead of considering a selection  of the form $u\int_{\O}k(y)|u(t,y)|^p\,dy$, we consider a selection of the form $u\mathcal{R}(u)$ with $\mathcal{R}: dom(\mathcal{R})\to \R^+$ a positive functional satisfying:
$ \exists p,q\ge 1 \text{ and } c_p,\alpha_p, R_p,C_q, \alpha_q, R_q$ \text{ positive constants  such that }, 
\begin{align*}
\mathcal{R}(u)>c_p\nlp{u}{p}{\O}^{p\alpha_p} \quad \text{when}\quad \nlp{u}{p}{\O}\ge R_p,  \\
\mathcal{R}(u)<C_q\nlp{u}{q}{\O}^{q\alpha_q} \quad \text{when}\quad \nlp{u}{q}{\O}\le R_q.  
\end{align*} 
A simple example of such $\mathcal{R}$ is the functional $\mathcal{R}(u):=\nlp{u}{p}{\O}^p\nlp{u}{q}{\O}^q$.

Similarly, the optimal existence criteria  Theorem \ref{msedp-thm2} will hold as well for a selection process  $u\mathcal{R}(x,u)$ such that 
$$\mathcal{R}_1(\cdot)\le \mathcal{R}(x,\cdot)\le \mathcal{R}_2(\cdot),$$
  where the $\mathcal{R}_i$ satisfy the above assumptions.
 
  We also wanted to stress that the regularity on the coefficient is far from optimal and extension  of our results for rougher coefficients $r,k,A$ should hold true. In order  to  keep our analysis of the asymptotic behaviour as simple as possible, we deliberately impose some regularity on the  considered coefficients. We believe that theses assumptions highlight the important point of the  method we used without altering the pertinence of the results obtained.

We also want to emphasize that these results are strongly related to the eigenvalue problem obtained by linearising the equation \eqref{msedp-eq-red} around the steady state $0$ which is a common feature  for classical reaction diffusion 
$$\frac{\partial u}{\partial t} =\Delta u +f(x,u),$$
where $f$ is a KPP type. However, the extension of Theorems \ref{msedp-thm1}, \ref{msedp-thm2} to unbounded domains $\O$ is far from obvious considering the multiplicity of notion of generalised eigenvalue \cite{BR}. Moreover, in  these situation the strict positivity of the kernel $k$ seems to introduce a strong dichotomy for the properties of the stationary solutions and consequently the dynamics observed for evolution problem. Indeed,  
in this direction some progress have recently been made for the so called nonlocal Fisher-KPP equation :
\begin{equation}
\frac{\partial u}{\partial t}=\Delta u + u(1-\phi\star u),
\label{msedp-eq-fishernoloc}
\end{equation}

where $\phi$ is a non-negative kernel. When $\phi$ is a positive integrable function, the constant $1$ is a positive solution. Moreover, for $\phi \in L^1\cap C^1$ positive so that $x^2\phi \in L^1$, it is shown in \cite{Berestycki2009}  that travelling semi-front exists for all speed $c\ge c^*$, i.e there exists $(U,c)$, so that $U>0$ and $U$ satisfies 
\begin{align*}
&U_{xx}+cU_x+U(1-\phi\star U)=0,\\
&\lim_{x\to +\infty} U=0,\quad \liminf_{x\to -\infty} U>0.
\end{align*} 
In particular when $c$ is large or $\phi$ is sufficiently concentrated or has a  positive Fourier transform, we have $\liminf_{x\to -\infty} U=\limsup_{x\to -\infty} U=1$, see \cite{Alfaro2012,Berestycki2009,Fang2011,Nadin2011}. On the contrary, from our analysis  the positive solution of  
  \begin{equation}\label{msedp-eq-ubdomain}
  \frac{\partial u}{\partial t}=\Delta u + u\left(1-\int_{\R^n}u(t,y)\,dy\right),
  \end{equation}
 converges uniformly to $0$, which is actually the only non-negative stationary solution. 

We mention also a recent related study \cite{Alfaro2012a} on a spatial demo-genetic model 
\begin{equation}\label{msedp-eq-sdg}
  \frac{\partial u}{\partial t} (t,x,y)=\Delta u(t,x,y) +  u\left(r(x-By)-\int_{\R}u(t,x,y')\,dy'\right),
  \end{equation}
which can be viewed as an extension of \eqref{msedp-eq-intro} where a spatial local adaptation is taken into account. The interplay between the space variable $x$ and the phenotypical trait variable $y$  corresponding to local adaptation is modelled through the growth term $r(x-By)$ which is a function taking its maximum at $0$. Generalisation of \eqref{msedp-eq-sdg} have been studied in \cite{Arnold2012,Prevost2004}

The extension of Theorems \ref{msedp-thm1}, \ref{msedp-thm2} and \ref{msedp-thm3} for mutation-selection equations involving a mutation kernel such as   
\begin{align}
&\frac{\partial u}{\partial t}=u\left(r(x)-\int_{\O}k_\eps(y)|u|^p(t,y)\,dy\right)+\int_{\O}M(x,y)[u(t,y)-u(t,x)]\,dy\quad \text{ in } \quad \R^+\times\O \label{msedp-eq-noloc}
\end{align}  
is still a work in progress. However, although the technique and tools developed in this article  are quite robust and can be applied in many situation, the lack of regularity of the positive solutions to \eqref{msedp-eq-noloc} introduces some strong difficulty that cannot be easily overcome.  Moreover, it has been proved by the author that such nonlocal problem can generates blow up phenomena, i.e.  $u(x,t)\rightharpoonup \delta_{x_0}+g$ with  $\delta_{x_0}$ the Dirac mass and $g$  a singular $L^1$ function. This blow up phenomena is in accordance with a recent result showing that in some situation the only stationary solution to \eqref{msedp-eq-noloc} are positive measure having a non-zero singular part \cite{Coville2013c}. The understanding of the long time behaviour of the positive solution to \eqref{msedp-eq-noloc} require then the development of new  analytical tools in order to analyse these blow-up phenomena.

 This paper is organised as follows. The Section \ref{msedp-section-entro} is dedicated  to  the nonlinear relative entropies and some functional inequalities that we will frequently use along this article. Next, we prove in Section \ref{msedp-section-blind} the Theorem \ref{msedp-thm1}. Finally in Section \ref{msedp-section-sta} and \ref{msedp-section-asb} we prove the existence of positive steady states (Theorem \ref{msedp-thm2}) and the global stability (Theorem \ref{msedp-thm3}). A construction of a smooth positive solution to the Cauchy problem is made in the appendix.

\section{Non-linear relative entropy identities and related functional inequality}\label{msedp-section-entro}
In this section we first establish a general identity which can be assimilated to a nonlinear relative entropy principle. We consider a parabolic  equation of the form
  
\begin{align}
&\frac{\partial u}{\partial t}(t,x)=u(t,x)(r(x)-\Psi(x,u)(t))+\nabla\cdot(A(x)\nabla u(t,x) )\quad \text{ in } \quad \R^+\times\O, \label{msedp-eq-gen}\\
&\frac{\partial u}{\partial n}(t,x)=0, \quad \text{ in } \quad \R^+\times\partial\O\label{msedp-eq-gen-bc}
\end{align}
where $\Psi(x,u)(t)$ denotes $\Psi(x,u)(t):=\int_{\O}K(x,y)|u|^p(t,y)\,dy.$
Then for any solution of \eqref{msedp-eq-gen}--\eqref{msedp-eq-gen-bc} we have  
\begin{theorem}[General Identity]\label{msedp-thm-gen-id}
Let $H$ be a smooth (at least $C^2$)  function. Let $\bar u >0$ and $u$  be two smooth solutions  of \eqref{msedp-eq-gen}--\eqref{msedp-eq-gen-bc}. Assume further that $\bar u$ is a stationary solution of \eqref{msedp-eq-gen}--\eqref{msedp-eq-gen-bc}.   Then we have 
 
\begin{equation}
\frac{d \oph{H}{\bar u}{u}(t)}{dt}=-\D(u)+\int_{\O} \bar u(x)H^{\prime}\left(\frac{u}{\bar u}(t,x)\right) \Gamma (t,x)u(t,x)\,dx
\end{equation}
 where $ \oph{H}{\bar u}{u}(t)$, $\D$ are the following quantity:
\begin{align*}
&\Gamma(t,x):=\Psi(x,\bar u) -\Psi(x,u)\\
&\oph{H}{\bar u}{u}(t):=\int_{\O}\bar u^2(x)H\left(\frac{u(x)}{\bar u(x) }\right)\, dx\\
&\D(u):=\int_{\O} \bar u^2(x)H^{\prime\prime}\left( \frac{u(x)}{\bar u(x)}\right) \left(\nabla\left(\frac{u}{\bar u}\right)\right)^t A(x)\nabla\left(\frac{u}{\bar u}\right) \, dx 
 \end{align*}
where $(\vec a)^t $ denotes the transpose of a vector of  $\R^N$. 
\end{theorem}

\dem{Proof:}
By \eqref{msedp-eq-gen}, by defining $\Gamma(t,x):=\Psi(x,\bar u(x)) -\Psi(x,u(t,x))$ we have 

 \begin{equation}
 \frac{\partial u}{\partial t}= \left(r(x) -\Psi(x,\bar u)u  +\nabla\cdot(A(x)\nabla u)\right) +\Gamma(t,x) u(x)
 \end{equation}
Using that $\bar u$ is also a stationary solution, we have for all $x$
$$(r(x)-\Psi(x,\bar u) )\bar u=-\nabla\cdot(A(x)\nabla \bar u),$$
and we can rewrite the above equation  as follows
  $$
  \frac{\partial  u(x)}{\partial t}=   \nabla\cdot(A(x)\nabla u) -\frac{u}{\bar u}\nabla\cdot(A(x)\nabla \bar u)              +\Gamma(t,x) u(x)
  $$
  By multiplying the above equality by $\bar u(x) H^{\prime}\left(\frac{u(x)}{\bar u(x)}\right)$ and by integrating over $\O$ we achieve 
\begin{multline}
\int_{\O} \bar u(x)H^{\prime}\left(\frac{u(x)}{\bar u(x)}\right) \frac{\partial u(x)}{\partial t}\, dx=\int_{\O} \bar u(x)H^{\prime}\left(\frac{u(x)}{\bar u(x)}\right) \Gamma(t,x) u(x)\, dx\\  +\int_{\O} H^{\prime}\left(\frac{u(x)}{\bar u(x)}\right) [\bar u(x)\, \nabla\cdot(A(x)\nabla u) -u(x)\, \nabla\cdot(A(x)\nabla \bar u(x))]\, dx .
\end{multline}
By integrating by part the last term and  rearranging the terms, it follows that 
\begin{multline}
\int_{\O} \bar u(x)H^{\prime}\left(\frac{u(x)}{\bar u(x)}\right) \frac{\partial u(x)}{\partial t}\, dx=\int_{\O} \bar u(x)H^{\prime}\left(\frac{u(x)}{\bar u(x)}\right) \Gamma(t,x) u(x)\, dx\\  -\int_{\O} \bar u^2(x)H^{\prime\prime}\left( \frac{u(x)}{\bar u(x)}\right) \left(\nabla\left(\frac{u}{\bar u}\right)\right)^tA(x) \nabla\left(\frac{u}{\bar u}\right) \, dx .
\end{multline}
Hence, we have  
$$\frac{d \oph{H}{\bar u}{u}(t)}{dt}=\int_{\O} \bar u(x)H^{\prime}\left(\frac{u(x)}{\bar u(x)}\right) \Gamma(t,x) u(x)\, dx-\D(u). $$
\fdem
\begin{remark}\label{msedp-rem-id}
We want to stress that if we replace $\bar u$ by any positive function $\tilde u$ satisfying
\begin{align*}
&\nabla\cdot(A(x)\nabla \tilde u(x) )=-\tilde u(x)\left(r(x)-\tilde \Psi(x,\tilde u)(t)\right)\quad \text{ in } \quad \O, \\
&\frac{\partial \tilde u}{\partial n}(x)=0, \quad \text{ in } \quad \partial\O
\end{align*}
it will affect the equality in Theorem \ref{msedp-thm-gen-id} only through the term $\Gamma$  which will be transform into
 $$\Gamma(t,x) =\tilde \Psi(x,\tilde u(x))-\Psi(x,u(t,x)).$$
\end{remark}
\begin{remark}
Under the  extra assumption $\frac{u}{\bar u} \in L^{\infty}(\O)$, we remark that the formulas will holds as well if we consider  homogeneous Dirichlet boundary conditions instead of Neumann boundary conditions. It is worth noticing that this extra condition is always satisfied in the Neumann case since for all positive stationary solution with homogeneous Neumann Boundary condition, we can show that  $\inf_{\bar \O}\bar u>0$.
\end{remark}

\begin{remark}
We remark that the above formula do not require any particular assumption on the $\Psi$ and as a consequence no particular assumption on the kernel $K$. Thus the formula holds as well for $K(x,y)=\delta_0$, which turns the equation \eqref{msedp-eq-gen} into a semi-linear PDE. In particular  when $\Psi(x,u)$ is  independent of $u$ $i.e p=0, K=\delta_0 $ then the formula in Theorem  \ref{msedp-thm-gen-id} is known as the standard relative entropy principle for linear equations see \cite{Michel2005}. 
\end{remark}

Next we establish a  useful functional inequality satisfied by vectors  $h \in \bar v^{\perp}$ where $\bar v^{\perp}$ denotes the linear subspace of $H^1(\O)$:
 
  $$\bar v^{\perp}:=\left\{h \in H^1(\O)\,\left |\,  \int_{\O}h\bar v =0, \quad \bar v\nabla h\cdot n-h\nabla \bar v\cdot n=0\quad \text{ on }\quad \partial \O \right. \right \} $$

\begin{lemma} \label{msedp-lem-fcineq}
Let $\bar v$ be a smooth ($C^{1,\alpha}(\O)$) positive bounded function in $\O$, so that $\inf_{\bar\O}\bar v>0$. Then there exists $\rho_1>0$ so that for all $h\in \bar v^{\perp}$  
$$\rho_1\nlp{h}{2}{\O}^2\le \int_{\O}\bar v^2\left(\nabla\left( \frac{h} {\bar v}\right)\right)^t A(x)\nabla\left( \frac{h} {\bar v}\right).$$

Moreover $\rho_1=\lambda_2$ where $\lambda_2$ is the second eigenvalue  of the linear eigenvalue problem
\begin{align*} 
&\nabla\cdot\left(A(x)\bar v^2\nabla\left( \frac{h} {\bar v}\right)\right)=-\lambda h\bar v \quad \text{ in }\quad  \O\\
&\bar v \frac{\partial h}{\partial n} -h\frac{\partial\bar v}{\partial n}=0 \quad \text{ in }\quad  \partial\O
\end{align*}
\end{lemma}
\dem{Proof :}
Let $\I$ be the following functional in $H^1(\O)$,
\begin{equation}
\I(h):=\frac{1}{\nlto{h}^2}\int_{\O}\bar v^2\left(\nabla\left( \frac{h} {\bar v}\right)\right)^t A(x)\nabla\left( \frac{h} {\bar v}\right). \label{msedp-eq-functional}
\end{equation}
Observe that from the homogeneity of the  $L^2$ norm we have 
\begin{equation}
\inf_{h\in \bar v^{\perp}, \nlto{h}=1}\I(h)=\inf_{h\in \bar v^{\perp}}\I(h),
\end{equation}
and the first part of the Lemma is proved if we show that   
\begin{equation}
\inf_{h\in \bar v^{\perp}, \nlto{h}=1}\I(h)>0,\label{msedp-eq-energy-expdecay}
\end{equation}

Let $d\mu$ denotes the positive measure $\bar v^2 dx$, then by construction $d\mu$ is absolutely  continuous with respect to the Lebesgue measure and vice versa.  So the Hilbert functional spaces $L^2_{d\mu}$ and $H^1_{d\mu}$ below are well defined :
\begin{align*}
&L^{2}_{d\mu}(\O):=\left\{u \, \left|\, \int_{\O}u^2(x)d\mu(x)<+\infty\right.\right\},\\
&H^1_{d\mu}(\O):=\left\{u\in L^{2}_{d\mu}(\O) \, \left|\, \int_{\O}|\nabla u|^2(x)d\mu(x)<+\infty\right.\right\}.
\end{align*}
Moreover the Rellich-Kondrakov compact embedding   $H^1_{d\mu}(\O)\hookrightarrow L^{2}_{d\mu}(\O)$ holds \cite{GT}. 
To obtain \eqref{msedp-eq-energy-expdecay}, we  argue as follows.  
 Let $(h_n)_{n\in \N}$ be a minimising sequence, by \eqref{msedp-eq-energy-expdecay} we can take $(h_n)_{n\in \N}$  so  that  $h_n\in\bar v^{\perp}, \nlto{h_n}=1$ for all $n$.  Let $g_n:=\frac{h_n}{\bar v}$, then by straightforward computation,
from \eqref{msedp-eq-functional} -- \eqref{msedp-eq-energy-expdecay}, we see that   $(g_n)_{n\in\N}$ is a minimising sequence of the functional
$$\j(g):=\frac{1}{\nlmu{g}{\O}}\int_{\O} (\nabla(g))^tA(x)\nabla(g)\,d\mu,$$
satisfying for all $n$, $ \nlmu{g_n}{\O}=\nlto{h_n}=1$. Moreover, we have for all $n$, $\frac{\partial g_n}{\partial n}=0$ on $\partial\O$ and
\begin{equation}
\int_{\O}g_n(x)\,d\mu(x)=\int_{\O}h_n(x)\bar v(x)\, dx=0.
\label{msedp-eq-constr1}
\end{equation}
We can also easily verify that 
$$ \inf_{h\in \bar v^{\perp}, \nlto{h}=1}\I(h)=\inf_{g\in H^1_{d\mu}, \int_{\O}g\,d\mu=0 }\j(g).  $$

 By construction the sequence $(g_n)_{n\in\N}$ is uniformly bounded in $H^1_{d\mu}(\O)$ and thanks to  Rellich-Kondrakov compact embedding,   there exists a subsequence $(g_{n_k})_{k\in N}$ which converges weakly in $H^1_{d\mu}(\O)$ and strongly in $L^2_{d\mu}(\O)$ to some $\tilde g \in H^1_{d\mu}(\O)$. Moreover, $\tilde g$ is a weak solution of 
\begin{align}
 &\nabla\cdot\left(A(x)\bar v^2\nabla\left( \tilde g\right)\right)=-\lambda \tilde g \bar v^2 \quad \text{ in }\quad  \O,\label{msedp-eq-el-lim}\\ 
 &\frac{\partial \tilde g}{\partial n}=0\label{msedp-eq-el-lim-bc} 
\end{align}
for some $\lambda \in \R$. Furthermore $\tilde g$ satisfies 
\begin{equation}
\int_{\O}\tilde g(x)\,d\mu(x)=0.
\label{msedp-eq-constr2}
\end{equation}

 Now  assume that $\lambda=0$, then the above equations \eqref{msedp-eq-el-lim}--\eqref{msedp-eq-constr2}  enforce  $\tilde g =0$ leading to the  contradiction $0=\nlmu{\tilde g}{\O}=1.$ Therefore $\lambda \neq 0$ and \eqref{msedp-eq-energy-expdecay} holds.
 
Now, since $A(x)$ and $\bar v$ are smooth and  $\mu$ is absolutely  continuous with respect to the Lebesgue measure, by standard elliptic regularity  we have $\tilde g \in   C^{2,\alpha}(\O)$ for some $\alpha$ and the function $\tilde h:=\bar v \tilde g \in C^{2}$ satisfies
\begin{align*} 
&\nabla\cdot\left(A(x)\bar v^2\nabla\left( \frac{\tilde h} {\bar v}\right)\right)=-\lambda \tilde h\bar v \quad \text{ in }\quad  \O,\\
&\int_{\O}\tilde h\bar v \,dx=0,\\
&\bar v \frac{\partial h}{\partial n} -h\frac{\partial\bar v}{\partial n}=0 \quad \text{ in }\quad  \partial\O.
\end{align*}

Now by dividing \eqref{msedp-eq-el-lim} by $\bar v^2$ we get the following eigenvalue problem
  \begin{align*} 
&\frac{1}{\bar v^2}\nabla\cdot\left(A(x)\bar v^2\nabla g\right)=-\lambda g \quad \text{ in }\quad  \O,\\
&\frac{\partial g}{\partial n}=0 \quad \text{ in }\quad  \partial\O,
\end{align*}

From standard Theory \cite{GT} there exists a sequence $\lambda_1<\lambda_2<\lambda_3<\ldots $ of eigenvalue of the above problem. Moreover there exists  an orthonormal basis $\{\psi_k\}_{k=1}^{\infty}$ of  $L^2$, so that $\psi_k$ satisfies
\begin{align*} 
&\frac{1}{\bar v^2}\nabla\cdot\left(A(x)\bar v^2\nabla \psi_k\right)=-\lambda_k  \psi_k \quad \text{ in }\quad  \O,\\
&\frac{\partial \psi_k}{\partial n}=0 \quad \text{ in }\quad  \partial\O.
\end{align*}
 By setting $\phi_k:=\frac{\psi_k}{\bar v}$, we can check that 

 \begin{align} 
&\nabla\cdot\left(A(x)\bar v^2\nabla\left( \frac{\phi_k} {\bar v}\right)\right)=-\lambda_k  \phi_k\bar v \quad \text{ in }\quad  \O,\label{msedp-eq-spec1}\\
&\bar v \frac{\partial \phi_k}{\partial n} -\phi_k\frac{\partial\bar v}{\partial n}=0 \quad \text{ in }\quad  \partial\O.\label{msedp-eq-spec2}
\end{align}

Here since $(0,\bar v)$ is a solution to \eqref{msedp-eq-spec1}--\eqref{msedp-eq-spec2}  and $\bar v>0$, we  see that $\phi_1=\bar v$ and $\lambda_1=0$.  So 
$$\inf_{h\in \bar v^{\perp}, \nlto{h}=1}\I(h) = \lambda_2,$$
since the $\lambda_i$ are ordered and $\phi_2\in \bar v^{\perp}$. 
\fdem

\section{The Blind competition case:}\label{msedp-section-blind}
In this section we analyse the  asymptotic behaviour of a positive smooth solution to \eqref{msedp-eq-intro}--\eqref{msedp-eq-intro-ci} when the competition  kernel $K(x,y)$ is  independent of $x$, i.e  $K(x,y)=k(y)$ with $k$ satisfying  \eqref{msedp-hyp1}. 
As we expressed in Theorem \ref{msedp-thm1} that we recall below, in this situation the problem \eqref{msedp-eq-red}--\eqref{msedp-eq-red-bc} has a unique positive stationary solution which attracts all the trajectories initiated from any nonnegative and non zero initial data. More precisely, we prove 
 
 \begin{theorem}
Assume $A,r,k$ satisfy \eqref{msedp-hyp1} and  $p\ge 1$.
Let $\lambda_1$ be the first eigenvalue  of    the problem 
\begin{align}
&\nabla\cdot(A(x)\nabla \phi(x) )+r(x)\phi(x)=-\lambda\phi(x) \quad \text{ in } \quad \O,\label{msedp-eq-pev}\\
&\frac{\partial \phi(x)}{\partial n}=0\quad \text{ on } \quad \partial\O,
\end{align}
then we have the following asymptotic behaviour for any positive smooth solution $u(t,x)$ to \eqref{msedp-eq-red}--\eqref{msedp-eq-red-bc}
\begin{itemize}
\item if $\lambda_1\ge 0$, there is no positive stationary solution and $u(t,x)\to 0$ as $t\to \infty$ 
\item if $\lambda_1< 0$, then $$u(t,x) \to \mu \phi_1 $$ where 
$\mu=\left(\frac{-\lambda_1}{\int_{\O}k(y)|\phi_1(y)|^p\,dy}\right)^{\frac{1}{p}}$ and $\phi_1$ is the positive eigenfunction associated to $\lambda_1$ normalized by $\nl{\phi_1}{\O}=1$. 
\end{itemize}
\end{theorem}

In the sequel of this section  to simplify the presentation we introduce the notation
$$\Psi(u):=\int_{\O}k(z)|u(y)|^p\,dy.$$ 
Before proving  the Theorem,  we start by establishing some useful  Lemmas.

  \begin{lemma}\label{msedp-lem-steady}Assume $\lambda_1<0$, then  there exists $\mu>0$ so that $\mu \phi_1$ is a  positive stationary solution  of \eqref{msedp-eq-red}. 
\end{lemma}
\dem{Proof:}
Let us normalised $\phi_1$ by $\nl{\phi_1}{\O}=1$. Then, by plugging $\mu \phi_1$ in \eqref{msedp-eq-red}, we end up finding $\mu$ so that 
$$\Psi(\mu \phi_1)=-\lambda_1. $$ 
Thus for $\mu=\left(\frac{-\lambda_1}{\int_{\O}k(y)|\phi_1(y)|^p\,dy}\right)^{\frac{1}{p}}$, $\mu \phi_1$ is a stationary solution of \eqref{msedp-eq-red}.

\fdem.

Next, we establish some  useful identities. Namely, we show 
\begin{lemma}\label{msedp-lem-liap}Let $q\ge 1$ and $H$ be the smooth convex function $H(s):\,s\mapsto s^q$. Let $\bar u$ be a positive stationary solution of \eqref{msedp-eq-red}-\eqref{msedp-eq-red-bc}, then   a  positive smooth solution $u(t,x)$ of \eqref{msedp-eq-red}--\eqref{msedp-eq-req-bc} satisfies
\begin{equation}
\frac{d \oph{q}{\bar u}{u}(t)}{dt}=- q(q-1)\int_{\O}\left(\frac{u(t,x)}{\bar u(x)} \right)^{q-2}\bar u^2 \left(\nabla\left( \frac{u(t,x)} {\bar u(x)}\right)\right)^tA(x)\nabla\left( \frac{u(t,x)} {\bar u(x)}\right)\,dx  + q(\Psi(\bar u)-\Psi(u))\oph{q}{\bar u}{u}(t). \label{msedp-eq-estim1}
\end{equation}
where $\oph{q}{\bar u}{u}(t):=\int_{\O}\bar u^2(x)\left(\frac{u(t,x)}{\bar u(x)} \right)^q\,dx$.
Furthermore, the functional  $\f(u):=\log\left(\frac{\oph{q}{\bar u}{u}(t)}{\left(\oph{1}{\bar u}{u}(t)\right)^{q}}\right)$ satisfies: 
\begin{equation} \frac{d}{dt}\f(u)=-\frac{q(q-1)}{\oph{q}{\bar u}{u}(t)} \int_{\O}\left(\frac{u(t,x)}{\bar u(x)} \right)^q\bar u^2 \left(\nabla\left( \frac{u(t,x)} {\bar u(x)}\right)\right)^tA(x)\nabla\left( \frac{u(t,x)} {\bar u(x)}\right)\,dx.\label{msedp-eq-liap}\end{equation}
\end{lemma}
\medskip

\begin{remark}
Note that in the particular case of $H(s)=s^2$, $\oph{2}{\bar u}{u}=\nlto{u}^2$. So  we get a Lyapunov functional involving the $L^2$ norm of $u$ instead of  a weighted $L^q$ norm of $u$. 
Indeed, we have $$ \frac{\partial }{\partial t}\left(\log\left(\frac{\nlto{u}^2}{\left(\oph{1}{\bar u}{u} \right)^2}\right)\right)=-\frac{2}{\nlto{u}^2}\int_{\O}\bar u^2 \left(\nabla\left( \frac{u(t,x)} {\bar u(x)}\right)\right)^tA(x)\nabla\left( \frac{u(t,x)} {\bar u(x)}\right)\, dx.$$
\end{remark}

\dem{Proof:}
The identity   \eqref{msedp-eq-estim1} is a straightforward consequence of Lemma \ref{msedp-thm-gen-id}. Indeed, for $H(s):=s^q$, by the Theorem \ref{msedp-thm-gen-id} we have:
$$ 
 \frac{d \oph{q}{\bar u}{u}(t)}{dt}=-\D(u)+\int_{\O} \bar u (x)H^{\prime}\left(\frac{u(t,x)}{\bar u(x)}\right) \Gamma (x)u(x)\,dx
$$
 where $\Gamma$, $\D$ are the following quantity:
\begin{align*}
&\Gamma(u(t)):=\Psi(\bar u) -\Psi(u)\\
&\D(u):=\int_{\O} H^{\prime\prime}\left( \frac{u(t,x)}{\bar u(x)}\right) \bar u^2(x) \left(\nabla\left( \frac{u(t,x)} {\bar u(x)}\right)\right)^tA(x)\nabla\left( \frac{u(t,x)} {\bar u(x)}\right) \, dx
\end{align*}
By observing that $\bar u(x) u(x)H^{\prime}\left(\frac{u(t,x)}{\bar u(x)}\right)=q\oph{q}{\bar u}{u}(t)$ and that $\Gamma$ is independent of $x$, we see that 
 $$\frac{d \oph{q}{\bar u}{u}(t)}{dt}=-\D(u)+q\Gamma \oph{q}{\bar u}{u}(t),$$
and the formula \eqref{msedp-eq-estim1} holds.

To obtain \eqref{msedp-eq-liap},  we observe that by taking $q=1$ in the  formula \eqref{msedp-eq-estim1} we get    
$$ \frac{d \oph{1}{\bar u}{u}(t)}{dt}=\Gamma\oph{1}{\bar u}{u}(t).$$
Since  $\oph{1}{\bar u}{u}(t)=\int_{\O}u(t,x)\bar u(x)\,dx>0$ for all times we see that
\begin{equation} \frac{d}{dt}\log(\oph{1}{\bar u}{u}(t))= (\Psi(\bar u)-\Psi(u)).\label{msedp-eq-cdif1}
\end{equation}
Similarly,  since  $\oph{q}{\bar u}{u}(t)>0$ for all times we have also 
\begin{multline}
\frac{d}{dt}\log(\oph{q}{\bar u}{u}(t))= -\frac{q(q-1)}{\oph{q}{\bar u}{u}(t)} \int_{\O}\left(\frac{u(t,x)}{\bar u(x)} \right)^{q-2}\bar u^2(x) \left(\nabla\left( \frac{u(t,x)} {\bar u(x)}\right)\right)^tA(x)\nabla\left( \frac{u(t,x)} {\bar u(x)}\right)\,dx \\  +q(\Psi(\bar u)-\Psi(u)).\label{msedp-eq-cdif2}
\end{multline}

By combining \eqref{msedp-eq-cdif1} and \eqref{msedp-eq-cdif2}  we end up with 
$$ \frac{d }{d t}\left(\log\left(\frac{\oph{q}{\bar u}{u}(t)}{\left(\oph{1}{\bar u}{u}(t)\right)^{q}}\right)\right)=-\frac{q(q-1)}{\oph{q}{\bar u}{u}(t)} \int_{\O}\left(\frac{u(t,x)}{\bar u(t,x)} \right)^{q-2}\bar u^2(x) \left(\nabla\left( \frac{u(t,x)} {\bar u(x)}\right)\right)^tA(x)\nabla\left( \frac{u(t,x)} {\bar u(x)}\right)\,dx. $$
\fdem 

As a straightforward application of this Lemma, we  deduce the following  \textit{a priori} estimates on the solution of \eqref{msedp-eq-red}--\eqref{msedp-eq-red-ci}. Namely, we have

\begin{lemma}\label{msedp-lem-esti}
 Let $u(t,x)\in C^1((0,+\infty), C^{2,\alpha}(\O))$ be a positive solution of \eqref{msedp-eq-red}-\eqref{msedp-eq-red-bc}    then  for all $q\ge 1$ there exists a positive constant $c_q(q,u(x,1))<C_q(q,u(x,1))$ so that for all $t \ge 1$ 
 $$c_q\le \nlp{u}{q}{\O}\le C_q. $$
\end{lemma}

\dem{Proof:}
Let us first show that for all $q\ge 1$ then there exists $C_q(q,u(x,1))$ so that for all $t \ge 1$ 
 \begin{equation}
 \nlp{u}{q}{\O}\le C_q. \label{msedp-eq-blind-esti1}
 \end{equation}
  First, let us  obtain an upper bound for $u$ when $q=1$. By Lemma \ref{msedp-lem-liap}, we have 
$$\frac{d\oph{1}{\mu \phi_1}{u}(t)}{d t } = (\Psi(\mu\phi_1)-\Psi(u))\oph{1}{\mu \phi_1}{u}(t), $$
where $\mu \phi_1$ is the stationary solution constructed in Lemma \ref{msedp-lem-steady}.
By using the definition of $\Psi$ and $\oph{1}{\mu\phi_1}{u}(t)$, and H\"older's inequality , we have for some $c_0>0$
   $$\frac{d\oph{1}{\phi_1}{u}(t)}{d t } \le \left[\lambda_1-c_0 \left(\int_{\O}|u(t,y)|\,dy\right)^p\right]\oph{1}{\phi_1}{u}(t). $$
Since $\nlp{u}{1}{\O}\sim \oph{1}{\phi_1}{u}(t)$, we get for some $\tilde c_0$ 
$$\frac{d\oph{1}{\phi_1}{u}(t)}{d t } \le \left[\lambda_1-\tilde c_0 \left(\oph{1}{\phi_1}{u}(t)\right)^p\right]\oph{1}{\phi_1}{u}(t). $$
So $\oph{1}{\phi_1}{u}(t)$ satisfies a logistic differential inequation, therefore there exists $C_1(u(x,1))>0$ so that for all $t\ge 1$, 
\begin{equation}
\oph{1}{\phi_1}{u}(t)\le C_1.\label{msedp-eq-blind-esti2}
\end{equation}

Now we can get an upper bounded for $u$ for all $q\ge 1$. Indeed, let us assume that $q> 1$ then by a straightforward application of the Lemma \ref{msedp-lem-liap} we have for all $q> 1$ and for all  $t\ge 1$,
$$\oph{q}{\mu\phi_1}{u}(t)\le \left(\oph{1}{\mu\phi_1}{u}(t)\right)^q \left(\frac{\oph{q}{\mu \phi_1}{u}(1)}{\left(\oph{1}{\mu \phi_1}{u}(1)\right)^{q}}\right). $$
By using the homogeneity of the norm $\oph{q}{\mu\phi_1}{u}$ and \eqref{msedp-eq-blind-esti2} we see that   for all $q> 1$ and for all  $t\ge 1$,  
$$\oph{q}{\phi_1}{u}(t)\le \left(\oph{1}{\phi_1}{u}(t)\right)^q \left(\frac{\oph{q}{\phi_1}{u}(1)}{\left(\oph{1}{\phi_1}{u}(1)\right)^{q}}\right)\le C_1^q\left(\frac{\oph{q}{\phi_1}{u}(1)}{\left(\oph{1}{\phi_1}{u}(1)\right)^{q}}\right).$$

Since for $q\ge 1$ $\nlp{u}{q}{\O}\sim \oph{q}{\phi_1}{u}$, \eqref{msedp-eq-blind-esti1} holds.

To prove the lower bound for $u$, by H\"older's inequality, it is enough to have a lower bound for $\nlp{u}{1}{\O}$.
Recall that $\oph{1}{\mu\phi_1}{u}(t)$ satisfies
$$\frac{d\oph{1}{\mu \phi_1}{u}(t)}{d t } = \left(\Psi(\mu\phi_1)-\int_{\O}k(y)|u(t,y)|^p\,dy\right)\oph{1}{\mu \phi_1}{u}(t). $$
Since \eqref{msedp-eq-blind-esti1} holds for all $q\ge 1$, by interpolation there exits positive constants $C,\alpha$ so that for all $t>1$
$\nlp{u}{p}{\O}\le C \nlp{u}{1}{\O}^{\alpha}$. Therefore $\oph{1}{\mu\phi_1}{u}(t)$ satisfies for all $t>1$
$$\frac{d\oph{1}{\mu \phi_1}{u}(t)}{d t } \ge \left(\Psi(\mu\phi_1)-C^p|k|_\infty\oph{1}{\mu\phi_1}{u}^{\alpha p}\right)\oph{1}{\mu \phi_1}{u}(t). $$
By using the logistic character of the above differential inequation, we deduce that $\oph{1}{\mu \phi_1}{u}(t)\ge c_1(u(x,1))$ for all $t>1$.
\fdem

We are now in position to prove the Theorem \ref{msedp-thm1}.

\dem{Proof of Theorem \ref{msedp-thm1}:}
Let $u(t,x)\in C^{1}((0,+\infty), C^{2,\alpha}(\O))$ be a positive solution of \eqref{msedp-eq-red}--\eqref{msedp-eq-red-bc}.
Assume first that $\lambda_1< 0$. Since $u>0$ then $u$ is a sub-solution of 
 \begin{align}
&\frac{\partial v(t,x)}{\partial t}=\nabla\cdot(A(x)\nabla v(t,x) )+ r(x)v(t,x)\quad \text{ in } \quad \R^+\times\O\label{msedp-eq-lin} \\
&\frac{\partial v(t,x)}{\partial n}=0 \quad \text{ in } \quad \R^+\times\partial\O\label{msedp-eq-lin-bc}\\
&v(x,0)=u(1,x)\quad \text{ in } \quad \O.\label{msedp-eq-lin-ci}
\end{align} 
Since $\lambda_1>0$ and $u(1,x)\in L^{\infty}$, for a large constant  $Ce^{\lambda_1 t}\phi_1(x)$ is then a super-solution of \eqref{msedp-eq-lin}-- \eqref{msedp-eq-lin-ci} and by the parabolic maximum principle we have $$u(x,t)\le  Ce^{\lambda_1 t}\phi_1(x) \to 0 \quad \text{ as }\quad t\to \infty.$$

Now let us  assume that $\lambda_1=0$. In this situation, by Lemma \ref{msedp-lem-liap} and using Remark \eqref{msedp-rem-id},
we observe that for all $q\ge 1$ we have, 
$$\frac{d \oph{q}{\phi_1}{u}(t)}{dt}=- q(q-1)\int_{\O}\left(\frac{u(t,x)}{\phi_1(x)} \right)^{q-2}\phi_1^2 \left(\nabla\left( \frac{u(t,x)} {\phi_1(x)}\right)\right)^tA(x)\nabla\left( \frac{u(t,x)} {\phi_1(x)}\right)\,dx  -q\Psi(u))\oph{q}{\phi_1}{u}(t).$$ 
Therefore, since $\Psi(u)$ is non-negative, we get  $\nlto{\nabla u} \to 0$ and for all $q\ge 1$ $\nlp{u}{q}{\O} \to 0$ as $t\to +\infty$. Since the coefficients of the parabolic equation are uniformly bounded, by a bootstrap argument using the Parabolic regularity, we get  $\|u\|_{\infty}\to 0$ as $t\to \infty$.
 

Lastly, we assume  $\lambda_1<0$ and let us denote $<,>$ the standard scalar product of $L^2(\O)$. Let  $\bar u$ be the stationary solution of \eqref{msedp-eq-red}-- \eqref{msedp-eq-red-bc} constructed in Lemma \ref {msedp-lem-steady} , i.e $\bar u:=\mu\phi_1$.
Since for all $t>0$, the solution $u(t,x)\in L^{\infty}$, then we can decompose $u$ the following way:
$$u(t,x):=\lambda(t) \bar u(x) +h(t,x)$$ with $h$  so that $<\phi_1,h>=0$.

Substituting $u$ by this decomposition in \eqref{msedp-eq-red} and using the equation satisfied by $\bar u$ it follows that 
\begin{align}
\lambda'(t)\bar u(x) +\frac{\partial h(t,x)}{\partial t} 
&=(\lambda_1 -\Psi(u(t)))\lambda(t)\bar u(x)+ (r(x)-\Psi(u))h(t,x) +\nabla\cdot(A(x)\nabla(h(t,x))). \label{msedp-eq-red-decomp1}
\end{align} 
By multiplying the above equation by $h$ and integrating over $\O$, it follows that 

$$<\frac{\partial h(t)}{\partial t},h>= <(r(x)-\Psi(u))h +\nabla\cdot(A(x)\nabla(h)),h>.$$
where we use that $h$ is orthogonal to $\bar u$.
Thus since $\oph{2}{\bar u}{h}(t):=\nl{h(t)}{\O}^2$, we have 
$$<\frac{\partial h}{\partial t},h>= \frac{1}{2}\frac{d\oph{2}{\bar u}{h}(t)}{dt}=< (r-\Psi(u))h +\nabla\cdot(A\nabla(h)),h>.$$

By following the computation developed for the proof of Theorem \ref{msedp-thm-gen-id} with $H(s)=s^2$, we see that 
\begin{equation}\label{msedp-eq-h}
\frac{d\oph{2}{\bar u}{h}(t)}{dt}= - \int_{\O}\bar u^2(x) \left(\nabla\left( \frac{h(t,x)} {\bar u(x)}\right)\right)^tA(x)\nabla\left( \frac{h(t,x)} {\bar u(x)}\right)  +(\lambda_1 -\Psi(u(t))\oph{2}{\bar u}{h}(t).
\end{equation}
Since $\oph{2}{\bar u}{h}(t)\ge 0$ for all times, let us analyse separately   the two  situations:  $\oph{2}{\bar u}{h}(t)>0$ for all times $t$ or  there exists $t_0\in \R$ so that $\oph{2}{\bar u}{h}(t_0)=0$. In the latter case,  from the above equation we see that  we must have $\oph{2}{\bar u}{h}(t)=0$ for all $t\ge t_0$ and so for all $t\ge t_0$,  we  must have $u(t)=\lambda(t)\bar u $ almost everywhere. Hence from \eqref{msedp-eq-red-decomp1} we are reduced to analyse the following ODE equation
$$\lambda'(t)=\lambda(t)(\lambda_1 -\tilde \Psi(\lambda(t))) $$
 where $\tilde \Psi $ is the increasing locally Lipschitz function defined by $\tilde \Psi(s):=s^p\int_{\O}k(y)\bar u(y)^p\,dx$. 

Note that since  by Lemma \ref{msedp-lem-esti} we have 
\begin{equation}
\lambda(t)<\bar u,\bar u>=<\bar u,u>=\oph{1}{\bar u}{u}(t) \le C_1, \label{msedp-eq-blind-esti-lambda}
\end{equation}
 we have $\lambda(t)\ge 0$ for all times $t$. 
The above ODE is of logistic type with non negative initial datum therefore by a standard argumentation we see that  $\lambda(t)$ converges to  $\bar \lambda>0$ where $\bar \lambda$ is the unique solution of  $\tilde \Psi (\bar \lambda )=\lambda_1$. By construction we have $\tilde \Psi(1)=\lambda_1$, so we deduce that $\bar \lambda=1$.
Hence, in this situation,   $u$ converges  pointwise to $\bar u$  as time goes to infinity.
 
In the other situation,  $\oph{2}{\bar u}{h}(t)>0$ for all $t$ and we claim that 
\begin{claim} \label{msedp-cla-energy}
 $\oph{2}{\bar u}{h}(t)\to 0$ as $t\to +\infty$.
 \end{claim}
 Assume the Claim holds true then we can conclude  the proof by arguing as follows.
From  the decomposition $u(t,x)=\lambda(t)\bar u(x) +h(t,x)$, we can express  the function  $\oph{1}{\bar u}{u}(t)$  by $\oph{1}{\bar u}{u}(t)=<u,\bar u>=\lambda(t)<\bar u ,\bar u>$. Therefore by using Theorem \ref{msedp-lem-liap} we deduce that 
\begin{equation}
\lambda'(t)=\left(\lambda_1 -\Psi[\lambda(t)\bar u(x) +h(t,x)]\right)\lambda(t). \label{msedp-eq-lambda}
\end{equation}
By using the definition of $\Psi$ and the binomial expansion it follows that $\lambda$ verifies the following ODE
\begin{align*}
\lambda'(t)&=(\lambda_1 -\tilde \Psi(\lambda(t)))\lambda(t)+ \lambda(t)(\Psi(\lambda(t)\bar u(x))-  \Psi(\lambda(t)\bar u +h(t)))\\
&=(\lambda_1 -\tilde \Psi(\lambda(t)))\lambda(t)+ \lambda(t) \left(\sum_{i=1}^{p}\binom{i}{p}\lambda^i(t)\int_{\O}\bar u^ih^{p-i}(t,x)\,dx \right),
\end{align*}
where $\binom{i}{p}$ denotes the binomial coefficient.
 Now by using  $\nlto{h(t)}^2=\oph{2}{\bar u}{h}(t)\to 0$ and Lemma \ref{msedp-lem-esti}, by interpolation  we deduce that   $\nlp{h(t)}{q}{\O} \to 0$ for all $q\ge 1$.  Therefore, since $\bar u \in L^{\infty}$ and by \eqref{msedp-eq-blind-esti-lambda}  $\lambda$ is bounded, we have 
 $$\lim_{t\to \infty}\left(\sum_{i=1}^{p}\binom{i}{p}\lambda^i(t)\int_{\O}\bar u^ih^{p-i}(t,x)\,dx \right)=0.$$
Thus $\lambda$ satisfies
$$
\lambda'(t)=(\lambda_1 -\tilde \Psi(\lambda(t)))\lambda(t)+ \lambda(t) o(1),
$$
and as above we can conclude that $\lambda(t)\to 1$ and $u$ converges  to $\bar u$ almost everywhere.

 \fdem

\dem{Proof of   Claim \ref{msedp-cla-energy}:} 
Since $\oph{2}{\bar u}{h}(t)>0$ for all $t$,  from \eqref{msedp-eq-h} and by following the proof of Lemma \ref{msedp-lem-liap} we see that 
\begin{equation}
\frac{d}{dt}\log\left[\frac{\oph{2}{\bar u}{h}(t)}{\left(\oph{1}{\bar u}{u}(t)\right)^2}\right]=  -\frac{1}{\oph{2}{\bar u}{h}(t)}\int_{\O}\bar u^2(x) \left(\nabla\left( \frac{h(t,x)} {\bar u(x)}\right)\right)^tA(x)\nabla\left( \frac{h(t,x)} {\bar u(x)}\right)\,dx. \label{msedp-cla-dF}
\end{equation} 
Thus the function $\tilde F:= \log\left[\frac{\oph{2}{\bar u}{h}(t)}{\left(\oph{1}{\bar u}{u}(t)\right)^2}\right]$ is a decreasing smooth function.

First we observe that the claim is proved if there exists a sequence $(t_n)_{n\in\N}$ going to infinity so that 
$ \oph{2}{\bar u}{h}(t_n)\to 0$. Indeed, assume such sequence exists and let  $(s_k)_{k\in \N}$ be a sequence going to $+\infty$.  Then  there exists $k_0$ and a subsequence $(t_{n_k})_{k\in \N}$ of  $(t_n)_{n\in\N}$ so that  for all $k\ge k_0$,  we have $s_k \ge t_{n_k}$. Therefore from the monotonicity of $\tilde F$ we have for all $k\ge k_0$
$$\log\left[\frac{\oph{2}{\bar u}{h}(s_k)}{\left(\oph{1}{\bar u}{u}(s_k)\right)^2}\right]\le \log\left[\frac{\oph{2}{\bar u}{h}(t_{n_k})}{\left(\oph{1}{\bar u}{u}(t_{n_k})\right)^2}\right].$$
By letting $k$ to infinity in the above inequality, we deduce that 
$$\lim_{k\to \infty}\log\left[\frac{\oph{2}{\bar u}{h}(s_k)}{\left(\oph{1}{\bar u}{u}(s_k)\right)^2}\right]=-\infty,$$
which implies that $\oph{2}{\bar u}{h}(s_k)\to 0, $ since by Lemma \ref{msedp-lem-esti}  $(\oph{1}{\bar u}{u}(t_k))_{k\in \N}$ is uniformly bounded.
The sequence $(s_k)_{k \in \N}$ being chosen arbitrarily this implies that $\oph{2}{\bar u}{h}(t) \to 0 $ as $t\to +\infty$.  
\medskip

Let us now prove that such sequence $(t_n)_{n\in \N}$ exists. Let us assume   by contradiction  that  $\inf_{t\in \R^+}\oph{2}{\bar u}{h}(t)>0$. 

From the monotonicity and the smoothness  of $\tilde F$  we deduce that there is  $c_0\in \R$ so that
  $$\tilde F(h(t))\to c_0 \quad \text{ and }\quad \frac{d}{dt}\tilde F(h(t))\to 0 \quad\text{ as }\quad t \to +\infty.$$

Thus by Lemma \ref{msedp-lem-esti} and \eqref{msedp-cla-dF} it follows that 
 
\begin{equation} \lim_{t\to \infty} \int_{\O}\bar u^2(x) \left(\nabla\left( \frac{h(t,x)} {\bar u(x)}\right)\right)^tA(x)\nabla\left( \frac{h(t,x)} {\bar u(x)}\right)\,dx=0.\label{msedp-eq-lim}
\end{equation}

Since for all $t$, $h(t)\in \bar u^{\perp}$, $\oph{2}{\bar u}{h}(t)=\nlto{h(t)}^2 $ and $\bar u =\mu\phi_1\in C^{2,\alpha}$ is strictly positive in $\bar \O$, by combining \eqref{msedp-eq-lim} and the Lemma \ref{msedp-lem-fcineq}
we get the contradiction 
$$ 0<\lim_{t\to \infty}\nlto{h(t)}^2\le \frac{1}{\rho_1} \lim_{t\to \infty} \int_{\O}\bar u^2(x) \left(\nabla\left( \frac{h(t,x)} {\bar u(x)}\right)\right)^tA(x)\nabla\left( \frac{h(t,x)} {\bar u(x)}\right)\,dx=0.
$$
\fdem

\section{The general competition case: Existence of positive stationary solution}\label{msedp-section-sta}
In this section we investigate the existence of a positive stationary solution of \eqref{msedp-eq-gen} and prove Theorem \ref{msedp-thm2}.
That is we look for positive solution of
 \begin{align}
&\nabla\cdot(A(x)\nabla v )+ v\left(r(x)-\Psi(x,v)\right)=0\quad \text{ in } \quad \O,\label{msedp-eq-gen-sta} \\
&\frac{\partial v}{\partial n}(x)=0 \quad \text{ in } \quad \partial\O,\label{msedp-eq-gen-sta-bc}
\end{align} 
where $\Psi(x,v)=\int_{\O}K(x,y)|v(y)|^p\,dy$.
First observe that when $\lambda_1\ge 0$, then there is no positive solution of \eqref{msedp-eq-gen-sta}--\eqref{msedp-eq-gen-sta-bc}.
Indeed,  by multiplying  by $\phi_1$ the equation  \eqref{msedp-eq-gen-sta} and integrating by parts it follows that 
$$0=-\lambda_1\int_{\O}v(x)\phi_1(x)\,dx -\int_{\O}\Psi(x,v)v(x)\phi_1(x)\,dx,$$
which implies  $\lambda_1\int_{\O}v(x)\phi_1(x)\,dx=\int_{\O}\Psi(x,v)v(x)\phi_1(x)\,dx=0$ since $\Psi(x,v),v$ and $\phi_1$ are non negative.
Thus $v=0$ almost everywhere since $\phi_1>0$. 

Let us then assume that $\lambda_1<0$.  Let $k>0$  so that  the operator $\nabla\cdot(A(x)\nabla  )+ r(x)-k$ with Neumann boundary condition  is invertible in $C^{0,\alpha}(\O)$ and a positive solution of \eqref{msedp-eq-gen-sta}--\eqref{msedp-eq-gen-sta-bc} is a positive fixed point of the map  $T$ 

$$\begin{array}{rcl}
T: C^{0,\alpha}(\O)&\to&C^{0,\alpha}(\O)\\
\\v&\mapsto&Tv:=(\nabla\cdot(A(x)\nabla  )+ r(x)-k)_n^{-1}[\Psi(x,v)v-kv]\end{array}.$$

To check that $T$ has a positive fixed point we use a degree argument. Let $x_0\in \O$ be fixed and let 
 $K^s(x,y)$ be defined by  
$$
K^s(x,y):=sK(x,y)+(1-s)K(x_0,y).
$$
Let us now consider the homotopy  $H \in C([0,1]\times C^{0,\alpha}(\O), C^{0,\alpha}(\O))$ defined by
$$\begin{array}{ccl}
H: [0,1]\times C^{0,\alpha}(\O)&\to&C^{0,\alpha}(\O)\\
(s,v)&\mapsto&H(s,v):=(\nabla\cdot(A(x)\nabla  )+ r(x) -k)_n^{-1} [\Psi_s(x,v)v-kv].\end{array},$$

where $\Psi_s(x,v):=\int_{\O}K^s(x,y)|v|^p(y)\,dy$.

One can see that $H(1,.)=T$ and $H(0,.)=T_0$ where $T_0$ corresponds to the map 
$$\begin{array}{rcl}
T_0: C^{0,\alpha}(\O)&\to&C^{0,\alpha}(\O)\\
\\v&\mapsto&T_0v:=(\Psi_0(v)-k)(\nabla\cdot(A(x)\nabla  )+ r(x)-k)_n^{-1}v.\end{array}$$

Note that there exists an unique positive fixed point to $T_0$ which can be constructed  as in Section \ref{msedp-section-blind}. 

Before computing the degree of $T_1$,  we obtain   some  \text{ a priori } estimates on the fixed point of the map $H(\cdot,\cdot)$. That is some estimates on the positive solution to the equation
$$
(\nabla\cdot(A(x)\nabla  )+ r(x)-k)_n^{-1} [\Psi_s(x,v)v-kv]=v 
$$
which rewrites:
\begin{align}
&\nabla\cdot(A(x)\nabla v)+ r(x)v = \Psi_s(x,v)v \label{msedp-eq-homotop}\\
&\partial_nv=0 \quad \text{ on } \partial\O \label{msedp-eq-homotop-bc}
\end{align} 

\begin{lemma}\label{msedp-lem-esti1}
Let $v$ be a continuous non negative solution of  \eqref{msedp-eq-homotop}-\eqref{msedp-eq-homotop-bc}. Then either $v\equiv 0$ or $v>0$ and there exists $\bar c_1$ and $\bar C_1$ independent of $s$  so that $$\bar c_1\le \int_{\O}|v|^p(x)\,dx\le \bar C_1. $$
\end{lemma}

\dem{Proof:}
The strict positivity of the solution $v$ is a straightforward consequence of the strong maximum principle. Therefore either $v \equiv 0$ or $v>0$.  So let us assume that $v>0$ and then by multiplying by $v$ the equation  \eqref{msedp-eq-homotop} and integrating by parts we see that 
$$
\int_{\O}r(x)v^2(x)\, dx-\int_{\O} \left(\nabla v(x)\right)^tA(x)\nabla v(x)\,dx=\int_{\O}\Psi^s(x,v)v(x)^2\,dx\ge K_{min}\int_{\O}|v(y)|^p\,dy\int_{\O}v^2(x)\,dx,$$
where $K_{min}:=\min_{x,y\in\bar\O\times\bar\O}K(x,y)$.
Therefore we get 
$$\frac{\|r\|_{\infty}}{K_{min}} \ge \int_{\O}|v(y)|^p\,dy.$$
We also  get 
$$
 \int_{\O}r(x)v^2(x)\, dx-\int_{\O} \left(\nabla v(x)\right)^tA(x)\nabla v(x)\,dx \le K_{max}\int_{\O}|v(y)|^p\,dy \int_{\O}v^2(x)\,dx
$$
with $K_{max}:=\max_{x,y\in\bar\O\times\bar\O}K(x,y)$ which leads to 
$$ \frac{\lambda_1}{K_{max}}\le \int_{\O}|v(y)|^p\,dy.$$

\fdem

We are now in position to prove the existence of a positive solution to the equation \eqref{msedp-eq-gen-sta} by means of the computation of the topological degree of $T-id$ on a well chosen set $\mathcal{O}\subset C^{0,\alpha}(\O)$.
Let us choose positive constants $c_2$ and $C_2$ so that $c_2<\bar c_1$ and $C_2>\bar C_1$ where $\bar c_1$ and $\bar C_1$ are the constants obtained in  Lemma \ref{msedp-lem-esti1}.  Let $\O$ be the following open set
$$\mathcal{O}:=\left\{v\in C^{0,\alpha}(\O), v\ge 0\,|\,   c_2\le \int_{\O} v^p(x)\,dx \le C_2 \right\}$$  
and let us compute $deg(T-Id, \mathcal{O},0)$.
By Lemma \ref{msedp-lem-esti1}  for all $s \in [0,1]$  $ H(s,v)-v\neq 0$ on $\partial \mathcal{O}$. Therefore using that  $H(.,.)$ is an homotopy,  since $T$ is a compact operator,   we conclude that  $deg(T-Id, \mathcal{O},0)=deg(H(1,.)-Id,\mathcal{O},0)=deg(H(0,.)-Id,\mathcal{O},0).$
By construction,  from Section \ref{msedp-section-blind}, one can check that $deg(H(0,.)-Id,\mathcal{O},0)\neq 0 $ since  the map $T_0$ has a unique positive non degenerated fixed point.
Thus $deg(T-Id, \mathcal{O},0)\neq 0$ which shows that $T$ has a fixed point in $\mathcal{O}$.
\fdem


\section{Stability of the dynamics, convergence to the equilibria}\label{msedp-section-asb}
In this section we prove Theorem \ref{msedp-thm3}. That is to say,  we analyse  the stability under some perturbation of the dynamics established for \eqref{msedp-eq-red}--\eqref{msedp-eq-red-bc} in Section \ref{msedp-section-blind}.  More precisely we investigate  the global dynamics  of solution of

 \begin{align}\label{msedp-eq-pert}
&\frac{\partial u(x,t)}{\partial t}=u\left[r(x)-\int_{\O}k_\eps(x,y)|u|^{p}(y)\,dy\right]+ \nabla\cdot\left(A(x)\nabla u(t,x)\right)\quad \text{ in }\quad \O\times \R^+,\\
&\frac{\partial u}{\partial n}(t,x)=0 \qquad \text{ in }\quad \partial\O\times \R^{+,*},\label{msedp-eq-pert-bc}\\
&u(x,0)=u_0(x)\ge 0,\label{msedp-eq-pert-ci}
\end{align}
where $p=1$ or $2$ and $k_\eps(x,y):=k_0(y)+\eps k_1(x,y)$ with $\eps$  a small parameter. 
To obtain the asymptotic behaviour in this case,  we follow the strategy  developed  in Section \ref{msedp-section-blind}. Namely, we start by showing some \textit{a priori} estimates on the solution $u(t,x)$,  then we analyse the convergence by means of some  differential inequalities. For convenience, we dedicate a subsection to each essential part of the proof.  
 \smallskip
 
\subsection{A priori estimate }~\\ 
We start by establishing some useful differential inequalities. Namely we show that 
\begin{lemma}\label{msedp-lem-diffineq}
Assume that $A,r, k_i$ satisfies \eqref{msedp-hyp1} and let $\phi_1$ be the positive eigenfunction associated to $\lambda_1(\nabla\cdot\left(A(x)\nabla \right))+r(x)$ with Neumann boundary condition.  Let $q\ge 1$ and $H$ be the smooth convex function $H(s):\,s\mapsto s^q$. Then there exists $\eps_0$ so that for all $\eps\le \eps_0$ and   for all positive solution $u\in C^{1}((0,\infty),C^{2,\alpha}(\O))$ of \eqref{msedp-eq-pert}--\eqref{msedp-eq-pert-bc}, we have for $t>0$ 
\begin{align*}
&\frac{d\oph{q}{\phi_1}{u}(t)}{dt}\le - \D_{q,\phi_1}[u](t) +   q(-\lambda_1- \alpha_{\eps,-}(u)) \oph{q}{\phi_1}{u}(t)\\
&\frac{d\oph{q}{\phi_1}{u}(t)}{dt}\ge  - \D_{q,\phi_1}[u](t) +  q(-\lambda_1-\alpha_{\eps,+}(u)) \oph{q}{\phi_1}{u}(t)
\end{align*}
where \begin{align*}
&\D_{q,\phi_1}[u](t):=q(q-1)\int_{\O}\left(\frac{u(t,x)}{\phi_1(x)} \right)^{q-2}\phi_1^2(x)\left(\nabla\left( \frac{u(t,x)} {\phi_1(x)}\right)\right)^tA(x)\nabla\left( \frac{u(t,x)} {\phi_1(x)}\right)\,dx\\
&\oph{q}{\phi_1}{u}:=\int_\O \left(\frac{u(t,x)}{\phi_1(x)}\right)^{q}\phi_1^2(x)\,dx\\
&\alpha_{\eps,\pm}(u):=\int_\O \left(k_0(y)\pm\eps\|k_1\|_{\infty}\right)|u(t,y)|^p\,dy
\end{align*}

\end{lemma}

\dem{Proof:}

Observe that since $u$ is positive, from \eqref{msedp-eq-pert} it follows that 
\begin{align*}
\frac{\partial u}{\partial t}(t,x)&\le [r(x) -\alpha_{-,\eps}(u)]u(t,x) +\nabla\cdot\left(A(x)\nabla u(t,x)\right) ,\\
\frac{\partial u}{\partial t}(t,x)&\ge [r(x) -\alpha_{+,\eps}(u)]u(t,x) +\nabla\cdot\left(A(x)\nabla u(t,x)\right).
\end{align*}

Let $\bar \o^+_{\eps}$ and $\bar \o^-_\eps$ be the stationary solutions of the corresponding equations with homogeneous Neumann boundary condition:  
\begin{align*}
\frac{\partial \o^-_\eps(t,x)}{\partial t}&= [r(x) -\alpha_{-,\eps}(\o^-_\eps)]\o^-_\eps(t,x) +\nabla\cdot\left(A(x)\nabla \o^-_\eps(t,x)\right) ,\\
\frac{\partial \o^+_\eps(t,x)}{\partial t}&=[r(x) -\alpha_{+,\eps}(\o^+_\eps)]\o^+_\eps(t,x) +\nabla\cdot\left(A(x)\nabla \o^+_\eps(t,x)\right).
\end{align*}

Let $\eps$ small enough, says $\eps\le \frac{k_{0,min}}{2\|k_1\|_{\infty}}$, then  by construction  $\bar \o_{\eps}^{\pm}$ exists and   we have $\bar \o^{\pm}_{\eps}=\mu_{\eps}^{\pm}\phi_1$.
Now by arguing as in the proof of Theorem \ref{msedp-thm-gen-id}, we obtain 
\begin{align*}
&\frac{d\h^-_{_{H,\bar \o^-_\eps}}[u](t)}{dt}\le - \D_{H,\bar\o^-_\eps}[u](t)  + q[-\lambda_1-\alpha_{\eps,-}(u)] \h^-_{_{H,\bar \o^-_\eps}}[u](t), \\
&\frac{d \h^+_{_{H,\bar \o^+_\eps}}[u](t)}{dt}\ge - \D_{H,\bar\o^+_\eps}[u](t)  + q[-\lambda_1-\alpha_{\eps,+}(u)]\h^+_{_{H,\bar \o^+_\eps}}[u](t).
\end{align*}
where \begin{align*}
&\h^{\pm}_{_{H,\bar \o^{\pm}_\eps}}[u](t):=\int_{\O}(\bar \o_{\eps}^{\pm})^2(x) H\left(\frac{u(t,x)}{\bar \o_{\eps}^{\pm}(x)} \right)\,dx,\\
&\D_{H,\bar \o^{\pm}_\eps}[u](t):=\int_{\O}H^{\prime\prime}\left(\frac{u(t,x)}{\bar \o^{\pm}_{\eps}(x)} \right)(\bar \o^{\pm}_{\eps}(x))^2\left(\nabla\left( \frac{u(t,x)} {\bar \o^{\pm}_{\eps}(x)}\right)\right)^tA(x)\nabla\left( \frac{u(t,x)} {\bar \o^{\pm}_{\eps}(x)}\right)\,dx.
\end{align*}

By using that $\bar \o^{\pm}_{\eps}=\mu_{\eps}^{\pm}\phi_1$, the definition of $H$ and the homogeneity of $\oph{q}{ \mu_\eps^{\pm}\phi_1}{u}$, we deduce  that  
\begin{align*}
&\frac{d\oph{q}{\phi_1}{u}(t)}{dt}\le -\D_{q,\phi_1}[u](t)  +   q[-\lambda_1- \alpha_{\eps,-}(u)] \oph{q}{\phi_1}{u}(t),\\
&\frac{d\oph{q}{\phi_1}{u}(t)}{dt}\ge  - \D_{q,\phi_1}[u](t)   +  q[-\lambda_1-\alpha_{\eps,+}(u)] \oph{q}{\phi_1}{u}(t).
\end{align*}

\fdem
\medskip

Next, we derive some \textit{ a priori } estimates for the solutions $u\in C^1((0,\infty),C^{2,\alpha}(\O)) $of \eqref{msedp-eq-pert}--\eqref{msedp-eq-pert-bc}.

\begin{lemma}\label{msedp-lem-estigen1}
Assume that $A,r, k_i$ satisfies \eqref{msedp-hyp1}.  Then there exists $\eps_1$ so that  we have : 
\begin{itemize}
\item[(i)]   For all $q'\ge 1$ there exists $\bar c_{q'}<\bar C_{q'}$ so that for all $\eps\le \eps_1$ and for all positive continuous stationary solution $\bar u_\eps$ to \eqref{msedp-eq-pert}--\eqref{msedp-eq-pert-bc} 
$$\bar c_{q'}\le \nlp{\bar u_\eps}{q'}{\O}<\bar C_{q'}.$$  
\item[(ii)] There exists $0<\bar c_\infty<\bar C_\infty$, so that for all $\eps \le \eps_1$ and for  all continuous stationary solution  $\bar u_\eps$ to \eqref{msedp-eq-pert}--\eqref{msedp-eq-pert-bc} 
$$\bar c_{\infty}\le \bar u_\eps \le \bar C_\infty.$$  

\item[(iii)] For all $1\le q'\le p$, there exists $ 0< C_{q'}$, so that for all $\eps\le \eps_1$ and for all  $u_\eps\in C^1((0,\infty),C^{2,\alpha}(\O)) $  positive solution to \eqref{msedp-eq-pert}--\eqref{msedp-eq-pert-bc} there exists  $\bar t $ so that for all  $t\ge \bar t$ 
 $$ \nlp{u_\eps(t)}{q'}{\O}\le  C_{q'}.$$
 \item[(iv)] For $p=1$ or $p=2$ there exists  a positive constant $c_{1}$ , so that for all $\eps\le \eps_1$ and for all  $u_\eps\in C^1((0,\infty),C^{2,\alpha}(\O)) $  positive solution to \eqref{msedp-eq-pert}--\eqref{msedp-eq-pert-bc} there exists  $\bar t $ so that for all  $t\ge \bar t$ $$ \nlp{u_\eps(t)}{1}{\O}\ge  c_{1}.$$
\end{itemize}

\end{lemma}

\dem{Proof:}
Let us first observe that $(ii)$ is a straightforward consequence of $(i)$ since $\bar u_\eps$  satisfies an elliptic equation with uniformly bounded continuous coefficient with respect to $\eps$ and $\bar u_\eps$. 
To prove $(i)$, we first show the estimates for $q'=p$.
First let us observe that  by replacing $u_\eps$ by $\bar u_\eps$  and taking $q=1$ in the formulas of  Lemma \ref{msedp-lem-diffineq}, we get for $\eps\le \eps_0$
\begin{align*}
&0\le [-\lambda_1- \alpha_{\eps,-}(\bar u_\eps)] \oph{1}{\phi_1}{\bar u_\eps},\\
&0\ge (-\lambda_1- \alpha_{\eps,+}(\bar u_\eps)) \oph{1}{\phi_1}{\bar u_\eps}  .
\end{align*}

From the latter inequalities,  by using the positivity of $\bar u_\eps$ and $\phi_1$ it follows that 
\begin{align*}
& -\lambda_1\ge \int_{\O}(k_0(y)-\sigma)\bar u_\eps^p(y) \, dy \ge \inf_{x\in\O}(k_0(x)-\sigma)\nlp{\bar u_\eps}{p}{\O}^p,\\
&-\lambda_1\le \int_{\O}(k_0(y)+\sigma)\bar u_\eps^p(y) \, dy \le \sup_{x\in\O}(k_0(x)+\sigma)\nlp{\bar u_\eps}{p}{\O}^p,
\end{align*}
where $\sigma:=\eps \|k_1\|_{\infty}$. Let $\kappa_0:=\frac{\inf_{x\in\O}k_0(x)}{2}$ and choose  $\eps$ small enough, says so that $\eps< \frac{\kappa_0}{\|k_1\|_{\infty}}=:\eps',$ we achieve for all $\eps \le \eps'$ and all stationary solution $\bar u_\eps$

\begin{equation}
\left(\frac{-\lambda_1}{\|k_0\|_\infty+\eps_1\|k_1\|_{\infty})}\right)^{\frac{1}{p}}=:\bar c_p\le  \nlp{\bar u_\eps}{p}{\O}\le \bar C_p:=\left(\frac{-\lambda_1}{\kappa_0}\right)^{\frac{1}{p}}.\label{msedp-eq-esti0-eps}
\end{equation}

 Now recall that $\bar u_\eps$ satisfies the elliptic equation 
 \begin{align*}
 &\nabla\cdot\left(A(x)\nabla \bar u_\eps(x)\right) +\left(r(x) -\int_{\O}k_\eps(x,y) \bar u_\eps^p(y)\, dy\right)\bar u_\eps(x)  =0\quad \text{ in }\quad \O, \\
 &\frac{\partial \bar u_\eps(x)}{\partial n}=0\quad \text{ in }\quad \partial\O. 
 \end{align*}
 From \eqref{msedp-eq-esti0-eps},   the coefficients of this linear equation  are uniformly bounded in $L^{\infty}$ with respect to $\eps \in [0,\eps']$.  So by using the elliptic regularity and Sobolev's embedding \cite{Brezis2010}, we can show that  for all $q\ge 1$ there exists $C>0$ so that  
 $$\|\bar u_\eps\|_{W^{2,q}(\O)}\le C, $$
 with $C$ independent of $\eps$ and $\bar u_\eps$. 
 Thus there exists $C_\infty>0$ independent of $\bar u_\eps$, so that 
\begin{equation} 
 \|\bar u_\eps\|_\infty\le C_{\infty}.\label{msedp-eq-esti00-eps}
 \end{equation}

To obtain the  desired uniform lower bound $\bar c_q$,  a standard interpolation  argument can be used \cite{Brezis2010} combining  \eqref{msedp-eq-esti0-eps} and \eqref{msedp-eq-esti00-eps}. 
\medskip

Let us now prove $(iii)$. Let $\kappa_1:=\|k_0\|_\infty+\eps_1\|k_1\|_{\infty}$  and $\kappa_0:=\frac{\inf_{x\in\O}k_0(x)}{2}$ then by Lemma \ref{msedp-lem-diffineq}, since $\eps \le \eps'$ we get for all $q\ge 1$ and all $t>0$

\begin{multline*}
\frac{d\oph{q}{\phi_1}{u_\eps}(t)}{dt}\le -q(q-1) \int_{\O}\left(\frac{u_\eps(t,x)}{\phi_1(x)} \right)^{q-2}\phi_1^2(x)\left(\nabla\left( \frac{u(t,x)} {\phi_1(x)}\right)\right)^tA(x)\nabla\left( \frac{u_\eps(t,x)} {\phi_1(x)}\right)\,dx  \\ +  q[-\lambda_1- \kappa_0\nlp{u_\eps}{p}{\O}^p] \oph{q}{\phi_1}{u_\eps}(t).
\end{multline*}
Since $\h_1(u)\sim \nlp{u}{1}{\O}$, by  H\"older's inequality and  by choosing $q=1$ in the above inequality, it follows that 
\begin{equation}
\frac{d\oph{1}{\phi_1}{u_\eps}(t)}{dt}\le [-\lambda_1- \tilde\kappa_0\oph{1}{\phi_1}{u_\eps}(t)^p] \oph{1}{\phi_1}{u_\eps}(t) \quad \text{ in } \quad  (0,\infty).
\end{equation}
Using the logistic character of the above equation, there exists $t_1$  so that  $\oph{1}{\phi_1}{u_\eps}(t)\le \frac{-2\lambda_1}{\tilde \kappa_0}$ for all $t\ge t_1$.
A similar argument can be done  for $q=p$, thus  $\oph{p}{\phi_1}{u_\eps}(t)\le C_p$ for all $t\ge t_p$ and by interpolation we get for all $1\le q\le p$
\begin{equation}
\nlp{u_\eps}{q}{\O}\le C_q \quad \text{ for all } \quad t\ge t':=\sup\{t_1,t_p\}. \label{msedp-eq-esti1-eps}
\end{equation}

 To obtain the lower bound $(iv)$, it is enough to get an uniform lower bound for $\oph{1}{\phi_1}{u_\eps}(t)$. By Lemma \ref{msedp-lem-diffineq} we have 
\begin{equation} 
 \frac{d\oph{1}{\phi_1}{u_\eps}(t)}{dt}\ge \left(-\lambda_1- \|k\|_{\infty}\int_{\O}u_\eps^p(y)\,dy\right) \oph{1}{\phi_1}{u_\eps}(t). \label{msedp-eq-pert-esti-u}
\end{equation}
\subsection*{Case 1: $p=1$}
In this situation,  since  $\oph{1}{\phi_1}{u_\eps}(t)\sim \nlp{u_\eps}{1}{\O}$,  we deduce that 
$$
\frac{d\oph{1}{\phi_1}{u_\eps}(t)}{dt}\ge (-\lambda_1- \kappa_1\oph{1}{\phi_1}{u_\eps}(t)) \oph{1}{\phi_1}{u_\eps}(t),
$$
for some $\kappa_1>0$. Hence, there exists $\bar t$ so that   $\oph{1}{\phi_1}{u_\eps}(t)\ge \frac{-\lambda_1}{2\kappa_1}$ for all $t>\bar t$.

\subsection*{Case 2: $p=2$}
In this situation, let us rewrite  $u_\eps(x,t):=\mu_\eps(t)\phi_1(x)+g_\eps(t,x)$ with $g(t,x) \perp \phi_1$ in $L^2(\O)$.
Equipped with this decomposition, we have 
\begin{align}
&\oph{1}{\phi_1}{u}(t)=\mu_\eps(t) \label{msedp-eq-pert-esti-u1}\\
&\nlto{u(t)}^2=\oph{2}{\phi_1}{u}(t)=\mu_\eps^2(t)+\nlto{g_\eps(t)}^2\label{msedp-eq-pert-esti-u2}\\
&\frac{d \oph{2}{\phi_1}{g_\eps}(t)}{dt}=\frac{d \oph{2}{\phi_1}{u_\eps}(t)}{dt}-2\mu_\eps(t)\mu_\eps'(t)\label{msedp-eq-pert-esti-u3}
\end{align}
So from \eqref{msedp-eq-pert-esti-u}, we get
\begin{equation} 
 \mu_\eps'(t)\ge \left(-\lambda_1- \|k_\eps\|_{\infty}\mu_\eps^2(t) - \|k_\eps\|_{\infty}\nlto{g_\eps(t)}^2\right) \mu_\eps(t). \label{msedp-eq-pert-esti-u4}
\end{equation}

Now  by combining \eqref{msedp-eq-pert-esti-u1}, \eqref{msedp-eq-pert-esti-u3} and Lemma \ref{msedp-lem-diffineq} we see that 
\begin{multline}
\frac{d \oph{2}{\phi_1}{g_\eps}(t)}{dt}\le - 2\int_{\O}\phi_1^2(x)\left(\nabla\left( \frac{g_\eps(t,x)} {\phi_1(x)}\right)\right)^tA(x)\nabla\left( \frac{g_\eps(t,x)} {\phi_1(x)}\right)\,dx  +  \frac{d \log{\mu_\eps^2(t)}}{dt} \oph{2}{\phi_1}{g_\eps}(t)\\ +2[\alpha_{\eps,+}(u_\eps)-\alpha_{\eps,-}(u_\eps)] \left(\mu_\eps^2(t)+\oph{2}{\phi_1}{g_\eps}(t)\right)
\label{msedp-eq-pert-esti1}
\end{multline}
By Lemma \ref{msedp-lem-fcineq} and using \eqref{msedp-eq-esti1-eps} it follows that  for $t\ge t'$
\begin{multline}
\frac{d \oph{2}{\phi_1}{g_\eps}(t)}{dt}-\frac{d \log{\mu_\eps^2(t)}}{dt} \oph{2}{\phi_1}{g_\eps}(t)\le - (2\rho_1(\phi_1)-4\eps\|k_1\|_\infty C_2)\oph{2}{\phi_1}{g_\eps}(t)     \\ +4\eps\|k_1\|_\infty C_2 C_1^2.
\label{msedp-eq-pert-esti2}
\end{multline}
Let $\Sigma:=\{t\ge t'\, |\, \oph{2}{\phi_1}{g_\eps}(t)>0\}$, then  we have for all  $t\in \Sigma$
\begin{equation}
\frac{d}{dt}\left(\log{\left[\frac{\oph{2}{\phi_1}{g_\eps}(t)}{\mu_\eps^2(t)}\right]}\right)\le - (2\rho_1(\phi_1)-4\eps\|k_1\|_\infty C_2) + \frac{4\eps\|k_1\|_\infty C_2 C_1^2}{\oph{2}{\phi_1}{g_\eps}(t)}.
\label{msedp-eq-pert-esti3}
\end{equation}
By choosing $\eps$ small enough, say $\eps \le \eps":=\min \left\{\eps',\frac{\rho_1(\phi_1)}{4C_2\|k_1\|_{\infty}}\right\}$, and  by letting $\delta:= 4\|k_1\|_\infty C_2 C_1^2$, by \eqref{msedp-eq-pert-esti3} we achieve for all $t\in \Sigma$
\begin{equation}
\frac{d}{dt}\left(\log{\left[\frac{\oph{2}{\phi_1}{g_\eps}(t)}{\mu_\eps^2(t)}\right]}\right)\le - \rho_1(\phi_1) + \frac{\eps\delta}{\oph{2}{\phi_1}{g_\eps}(t)}.
\label{msedp-eq-pert-esti4}
\end{equation}
To obtain the lower bound, the proof follows now  three  steps:
\subsection*{Step One}
We claim that 
\begin{claim}\label{msedp-cla-esti-s1}
For all $\eps\le \eps^{''}$, there exists $t_0>t'$ so that 
$$\oph{2}{\phi_1}{g_\eps}(t_0)<\frac{2\delta \eps}{\rho_1(\phi_1)}.$$ 
\end{claim}
\dem{Proof:}
 
Assume by contradiction that for all $t\ge t'$ we have 
$$\oph{2}{\phi_1}{g_\eps}(t_0)\ge \frac{2\delta \eps}{\rho_1(\phi_1)}.$$ 

Therefore it follows from \eqref{msedp-eq-pert-esti4} that for all $t>t'$  
\begin{equation}
\frac{d}{dt}\left(\log{\left[\frac{\oph{2}{\phi_1}{g_\eps}(t)}{\mu_\eps^2(t)}\right]}\right)\le - \frac{\rho_1(\phi_1)}{2}. 
\label{msedp-eq-pert-esti5}
\end{equation}

  Thus $ F(t):=\log{\left[\frac{\oph{2}{\phi_1}{g_\eps}(t)}{\mu_\eps^2(t)}\right]}$ is a decreasing function which by assumption is bounded from below for all $t\ge t'$. Therefore $F$ converges as $t$ tends to $+\infty$ and $\frac{dF}{dt} \to 0$. 
  Hence for $t$ large enough, we get the contradiction 
 $$
-\frac{\rho_1(\phi_1)}{4}\le \frac{d}{dt}\left(\log{\left[\frac{\oph{2}{\phi_1}{g_\eps}(t)}{\mu_\eps^2(t)}\right]}\right)\le - \frac{\rho_1(\phi_1)}{2}.
$$
\fdem
\subsection*{Step Two}
Let  $\eps_1$ and $\gamma(t_0)$ be the following quantities 
\begin{align*}
&\eps_1:=\min\left\{\eps", \frac{-\lambda_1\rho_1(\phi_1)}{8\|k_\eps\|_\infty\delta}\right\},\\
&\gamma(t_0):=\min\left\{\mu_\eps(t_0),\sqrt{\frac{-\lambda_1}{2\|k_\eps\|_\infty}}\right\}
\end{align*}
and let $Q$ be the real map
$$
\begin{array}{rcl}
\R^+&\to&\R^+\\
x&\mapsto&A\frac{Bx}{Bx+C}
\end{array}
$$
 where $A:=\frac{-\lambda_1}{2\|k\|_{\infty}}, B:=\rho_1(\phi_1)$ and $C:=2\epsilon\delta$.
We claim that 
\begin{claim} \label{msedp-cla-esti-s2}
For all $\eps\le \eps_1$ we have 
\begin{itemize}
\item[(i)] For all $t\ge t_0$, 
$$\mu_\eps^2(t)\ge \gamma^2(t_0).$$
\item[(ii)] There exists $t^{\prime}_1\ge t_0$ so that for all $t>t^{\prime}_1$
$$\mu_\eps^2(t)\ge Q(\gamma^2(t_0)).$$
\end{itemize}
\end{claim}

\dem{Proof:}

Let us denote $\Sigma^{\pm}$ and $\Sigma_0$ the following sets
\begin{align*}
&\Sigma^+:=\left\{t \ge t_0\, |\, \oph{2}{\phi_1}{g_\eps}(t)>\frac{2\delta \eps}{\rho_1(\phi_1)}\right\},\\
&\Sigma^-:=\left\{t \ge t_0\, |\, \oph{2}{\phi_1}{g_\eps}(t)\le \frac{2\delta\eps}{\rho_1(\phi_1)}\right\},\\
&\Sigma_0:=\left\{t\ge t_0 \, |\, \mu_\eps(t)\ge \min\left\{\mu_\eps(t_0),\sqrt{\frac{-\lambda_1}{2\|k_\eps\|_\infty}}\right\}\right\}.\\
\end{align*}

By construction $[t_0,+\infty)=\Sigma^+\cup \Sigma^-, t_0\in \Sigma^-$ 
and for all $\eps \le \eps_1$ we have
$$-\lambda_1-\|k_\eps\|_\infty\frac{2\eps\delta}{\rho_1(\phi_1)}\ge -\frac{-\lambda_1}{2}. $$
Let us now prove $(i)$. Let $\tilde t_0$ be the following time
$$\tilde t_0:=\sup\{t\ge t_0 \,|\, [t_0,t]\subset \Sigma^-\}.$$
By continuity of $\oph{2}{\phi_1}{g_\eps}(t)$, it follows from $\oph{2}{\phi_1}{g_\eps}(t_0)<\frac{2\delta \eps}{\rho_1(\phi_1)}$  that  $\tilde t_0>t_0$. Moreover we deduce   from \eqref{msedp-eq-pert-esti-u4}   that $\mu_\eps$ satisfies on $(t_0,\tilde t_0)$:
\begin{equation} 
 \mu_\eps'(t)\ge \left(-\frac{\lambda_1}{2}- \|k_\eps\|_{\infty}\mu_\eps^2(t) \right) \mu_\eps(t). \label{msedp-eq-pert-esti-u5}
\end{equation}
Therefore $\mu_\eps(t)\ge \min\left\{\mu_\eps(t_0),\sqrt{\frac{-\lambda_1}{2\|k_\eps\|_\infty}}\right\}$ for $t\in [t_0,\tilde t_0)$ which enforces $(t_0,\tilde t_0)\subset \Sigma_0$. 
Let $t^*$ be the following quantity
$$t^*:=\sup\{t\ge t_0 \,|\, (t_0,t)\subset \Sigma_0\}.$$	
From above  $(t_0,\tilde t_1)\subset \Sigma_0$, so we have $t^*\in (t_0,+\infty]$.
We will show that  $t^*=+\infty$. If not,   $t^*<+\infty$ and from the above arguments we can see that 
$\oph{2}{\phi_1}{g_\eps}(t^*)\ge\frac{2\delta \eps}{\rho_1(\phi_1)}$.
By definition of $t^*$, we have the following dichotomy since $[t_0,+\infty)=\Sigma^+\cup \Sigma^-$:
\begin{itemize}  
\item $t^*\in \Sigma^-$ and there exists $t^*<t^{*,+}\in \Sigma^+$ so that $(t^*,t^{*,+})\subset \Sigma^+$
\item $t^*\in \Sigma^+$ and there exists $t^{*,-}<t^*<t^{*,+}$ so that $t^{*,-}\in\Sigma_0\cap \Sigma^-, t^{*,+} \in \Sigma^+$ and  $(t^{*,-},t^{*,+}]\subset \Sigma^+$
 \end{itemize}
 In both cases   we see from \eqref{msedp-eq-pert-esti5} that on  $(t^{*,-},t^{*,+}]$ the function $F(t)= \log{\left[\frac{\oph{2}{\phi_1}{g_\eps}(t)}{\mu_\eps^2(t)}\right]}$ is decreasing and we have for all $t\in (t^{*,-},t^{*,+}] F(t)<F(t^{*,-})$ which leads to 
 $$ \mu_\eps^2(t^{*,-})\le \mu_\eps^2(t)\frac{\oph{2}{\phi_1}{g_\eps}(t^{*,-})}{\oph{2}{\phi_1}{g_\eps}(t)}. $$
Thus we get for all $t\in (t^{*,-},t^{*,+}]$ $$ \gamma(t_0)\le \mu_\eps(t),$$
 since $t^{*,-}\in \Sigma^{-}\cap \Sigma_0$ and $t\in \Sigma^+$.
 As a consequence we have $t^*< t^{*,+}\in \Sigma_0,$ which contradicts the definition of $t^*$.

Hence  $t^*=\infty$ and 
\begin{equation}
\mu_\eps(t)\ge \min\left\{ \mu_\eps(t_0), \sqrt{\frac{-\lambda_1}{2\|k_\eps\|_\infty}}\right\} \quad \text{ for all }\quad t\ge t_0. \label{msedp-eq-pert-esti6}
\end{equation}

 Let us now prove $(ii)$. By arguing on each connected component of $\Sigma^+$,  since by \eqref{msedp-eq-pert-esti4}  $F(t)= \log{\left[\frac{\oph{2}{\phi_1}{g_\eps}(t)}{\mu_\eps^2(t)}\right]}$ is a decreasing function  one has for all $t\in \Sigma^+$
$$ \oph{2}{\phi_1}{g_\eps}(t)\le \frac{\mu_\eps^2(t)}{\gamma^2(t_0)}\frac{2\eps\delta}{\rho_1(\phi_1)}.$$ 
 By construction, from  \eqref{msedp-eq-pert-esti6} we also have  for all $t\in \Sigma^-$ 
 $$ \oph{2}{\phi_1}{g_\eps}(t)\le \frac{\mu_\eps^2(t)}{\gamma^2(t_0)}\frac{2\eps\delta}{\rho_1(\phi_1)}.$$ 
 Therefore for all $t\ge t_0$ we get 
 \begin{equation}
 \oph{2}{\phi_1}{g_\eps}(t)\le \frac{\mu_\eps^2(t)}{\gamma^2(t_0)}\frac{2\eps\delta}{\rho_1(\phi_1)}. \label{msedp-eq-pert-esti7}
\end{equation}
Now by combining \eqref{msedp-eq-pert-esti7} with \eqref{msedp-eq-pert-esti-u4} it follows that for all $t\ge t_0$, $\mu_\eps(t)$ satisfies 
$$ \mu_\eps^{\prime}(t)\ge \left[-\lambda_1 -\|k\|_{\infty}\mu_\eps^2(t)\left(1+\frac{2\eps\delta}{\rho_1(\phi_1)\gamma^2(t_0)}\right)\right] \mu_\eps(t).$$
Hence, by using the logistic character of the above equation we have for some  $t^{\prime}_1$ for all $t\ge t^{\prime}_1$
\begin{equation} \label{msedp-eq-pert-esti8}
\mu_\eps^2(t)\ge \frac{-\lambda_1}{2\|k\|_{\infty}}\frac{\rho_1(\phi_1)\gamma^2(t_0)}{\gamma^2(t_0)\rho_1(\phi_1)+2\eps\delta}=Q(\gamma^2(t_0)). 
\end{equation}
\fdem
\subsection*{Step Three}
Finally we claim that 
\begin{claim}
There exists $\bar t$ so that for all $t\ge \bar t$ 
$$\mu_\eps^2(t)\ge \frac{-\lambda_1}{8\|k\|_\infty}.$$
\end{claim}

\dem{Proof:}

By an elementary analysis, one can check that  the map $Q(x)=A\frac{Bx}{Bx+C}$ is monotone increasing and  has a unique positive fixed point $x_0=\frac{AB-C}{B}=\frac{-\lambda_1}{2\|k\|_{\infty}} -\frac{2\epsilon\delta}{\rho_1(\phi_1)}\ge \frac{-\lambda_1}{4\|k\|_{\infty}}>0$.  We can also check that  the iterated  map $Q^{n+1}(x):=Q(Q^n(x))$ satisfies for any $x^*\in (0,+\infty)$
\begin{equation}\label{msedp-eq-pert-Q} 
\lim_{n\to \infty}Q^{n}(x^*)=x_0.
\end{equation}
Now recall that by the previous step, we have for all $t\ge t^{\prime}_1$, 
$$\mu_\eps^2(t)\ge Q(\gamma^2(t_0))=Q\left(\min\left\{\mu_\eps^2(t_0), \frac{-\lambda_1}{2\|k\|_{\infty}}\right\}\right).$$
Since $Q$ is monotone increasing and $\frac{-\lambda_1}{2\|k\|_{\infty}}>x_0$ we deduce from \eqref{msedp-eq-pert-esti8} that for all $t\ge t^{\prime}_1$
\begin{equation}\label{msedp-eq-pert-esti9}
\mu_\eps^2(t)\ge \min\left\{x_0, Q(\mu_\eps^2(t_0))\right\}.
\end{equation}
By using now step one with $t^{\prime}_1$ instead of $t'$, it follows that there exists $t_1\ge t^{\prime}_1$ so that 
 $\oph{2}{\phi_1}{g_\eps}(t_1)< \frac{2\delta\eps}{\rho_1(\phi_1)}$.
We can then replace $t_0$ by $t_1$ in Step two, to obtain the existence of  $t^{\prime}_2>t_1$ so that for all $t\ge t^{\prime}_2$ we have  
$$\mu_\eps^{2}(t)\ge Q(\gamma(t_1)^2)=Q\left(\min\left\{\mu_\eps^2(t_1), \frac{-\lambda_1}{2\|k\|_{\infty}}\right\}\right),$$
which by using the monotonicity of $Q$, $\frac{-\lambda_1}{2\|k\|_{\infty}}>x_0$ and \eqref{msedp-eq-pert-esti9} leads to 
\begin{equation}\label{msedp-eq-pert-esti10}
\mu_\eps^2(t)\ge \min\left\{x_0, Q\left[\min\left\{x_0, Q(\mu_\eps^2(t_0)\right\}
\right]\right\}
\end{equation}
for all $t\ge t^{\prime}_2$.

Since $x_0$ is a fixed point of $Q$, it follows from \eqref{msedp-eq-pert-esti10} that  for all $t\ge t^{\prime}_2$

\begin{equation}\label{msedp-eq-pert-esti11}
\mu_\eps^2(t)\ge \min\left\{x_0, Q[Q(\mu_\eps^2(t_0))]\right\}= \min\left\{x_0, Q^2(\mu_\eps^2(t_0))\right\}.
\end{equation}

By arguing inductively, we can then construct an increasing sequence $(t^{\prime}_n)_{n\in \N_{0}}$ so that for all $n$ and for all $t\ge t^{\prime}_n$ we have 
\begin{equation}\label{msedp-eq-pert-esti12}
\mu_\eps^2(t)\ge  \min\left\{x_0, Q^n(\mu_\eps^2(t_0))\right\}.
\end{equation}
Since $\mu_\eps^2(t_0)>0$, by \eqref{msedp-eq-pert-Q} there exists $n_0$ so that $Q^n(\mu_\eps^0(t_0))\ge \frac{x_0}{4}=\frac{-\lambda_1}{8\|k\|_{\infty}}$. Hence, by \eqref{msedp-eq-pert-esti12}
we have for all $t\ge t^{\prime}_{n_0}$
$$ \mu_\eps^2(t)\ge \frac{-\lambda_1}{8\|k\|_{\infty}} .$$

\fdem
\medskip

Finally, we establish an estimate on $\rho_1(\bar u_\eps)$ where $\rho_1(\bar u_\eps)$ is the constant defined in Lemma \ref{msedp-lem-fcineq} for the positive vector $\bar u_\eps$.  
Namely, we show that 
\begin{lemma}\label{msedp-lem-esti-fci}
There exists $\bar \rho>0$, so that for all $\eps\in [0,\eps_1)$ and for all positive stationary solution $\bar u_\eps$ of \eqref{msedp-eq-pert}--\eqref{msedp-eq-pert-bc}, we have 
$$\rho(\bar u_\eps)\ge \bar \rho$$
\end{lemma}

\dem{Proof:}
From the proof of Lemma \ref{msedp-lem-fcineq},  if we let  $d\mu_\eps$, $L^2_{\mu_{\eps}}$  and $H^1_{\mu_{\eps}}$ be  respectively the positive measure $d\mu_{\eps}=\bar u_\eps^2 dx$,  the following functional space: 
\begin{align*}
&L^{2}_{\mu_{\eps}}(\O):=\left\{u \, \left|\, \int_{\O}u^2(x)d\mu_\eps(x)<+\infty\right.\right\}\\
&H^1_{\mu_{\eps}}(\O):=\left\{u\in L^{2}_{\mu_{\eps}} \, \left|\, \int_{\O}|\nabla u|^2(x)d\mu_\eps(x)<+\infty\right.\right\}
\end{align*}
we have 
$$ 0<\rho(\bar u_\eps)=\inf_{g\in H^1_{d\mu_{\eps}}, \int_{\O}g\,d\mu_{\eps}=0 }\j(g)  $$
with $\j$ the functional
$$\j(g):=\frac{1}{\|g\|_{L^{2}_{\mu_{\eps}}(\O)}}\int_{\O}\left(\nabla(g)\right)^tA(x)\nabla(g)\,d\mu_\eps.$$
Let $$\nu:=\inf_{d\mu_{\eps}=\bar u_\eps^2 dx }\rho(\bar u_\eps),$$ where $\eps\in [0,\eps_1]$ and  $\bar u_\eps$  is any stationary solution of \eqref{msedp-eq-pert}--\eqref{msedp-eq-pert-bc}, then we have 
$$\rho(\bar u_\eps)\ge \nu \ge 0.$$
We claim that  $\nu>0$. Indeed, if not then there exists a sequence of positive measure $\bar u_n^2 dx$ so that 
 $$ \lim_{n\to \infty}\rho(\bar u_n)=0.$$
 Since $0\le \eps \le \eps_1$,  by Lemma \ref{msedp-lem-estigen1} the sequence $(\bar u_n)_{n\in\N}$ is uniformly bounded in $W^{2,q}(\O)$ for all $q\ge 1$. Therefore by the Rellich-Kondrakov Theorem, there exists  a subsequence $(\bar u_{n_k})_{k\in\N}$ which converges to $\tilde u$ a non-negative  solution of  \eqref{msedp-eq-pert}--\eqref{msedp-eq-pert-bc} for some $\bar \eps$. By Lemma  \ref{msedp-lem-estigen1}, we see also that $\tilde u$ is non trivial and positive.
 Thus by applying Lemma \ref{msedp-lem-fcineq} with $\tilde u$ we get the contradiction 
 $$ 0<\rho(\tilde u)=0.$$ 
 \fdem
\subsection{Asymptotic Behaviour}~\\
We are now in position to obtain the asymptotic behaviour of the solution $u_\eps(t,x)$ as $t$ goes to $+\infty$ for  $\eps \in [0,\eps^*]$, where $\eps^*$ is to be determined later on. 

Let us first introduce some practical notation:
\begin{align*}
&\Psi_0(v):=\int_{\O}k_0(y)|v(y)|^p\,dy, \;  \Psi_1(x,v):=\int_{\O}k_1(x,y)|v(y)|^p\,dy,   \;        \Psi_{\eps}(x,v):=\Psi_0(v)+\eps\Psi_1(x,v)\\ 
&\tilde \Psi_{\eps}(v):=\int_{\O}\Psi_\eps(x,v)v^2(x)\, dx.
\end{align*}

When  $\lambda_1\le 0$, then the proof of  Section \ref{msedp-section-blind} holds as well for solution of \eqref{msedp-eq-pert} --  \eqref{msedp-eq-pert-ci} and $u(t,x)\to 0$ as $t\to 0$. 
So let us  assume  $\lambda_1<0$ and let us denote $<,>$ the standard scalar product of $L^2(\O)$. Let  $\bar u_\eps$ be a positive stationary solution of \eqref{msedp-eq-pert}-- \eqref{msedp-eq-pert-bc}. Such solution exists from  Section \ref{msedp-section-sta}.
Since for all $t>0$ the solution $u_\eps(t,x)\in L^{2}$, we can decompose $u_\eps$ as follows:
$$u_\eps(t,x):=\lambda_\eps(t) \bar u_\eps +h_\eps(t,x)$$ with $h_\eps$  so that $<\bar u_\eps,h_\eps>=0$.

From this decomposition and  by using Theorem \ref{msedp-thm-gen-id} we get:  
\begin{align}
&\lambda_\eps(t) <\bar u_\eps,\bar u_\eps >=\oph{1}{\bar u_\eps}{u_\eps}(t) ,\label{msedp-eq-sysgen-asb1}\\
&\frac{d\oph{2}{\bar u_\eps}{h_\eps}(t)}{dt}=\frac{d\oph{2}{\bar u_\eps}{u_\eps}(t)}{dt} -2\lambda\lambda' <\bar u_\eps,\bar u_\eps>\label{msedp-eq-sysgen-asb2}\\
&\lambda^{\prime}_\eps(t) <\bar u_\eps,\bar u_\eps >=\int_{\O}(\Psi_\eps(x,\bar u_\eps)-\Psi_\eps(x,u_\eps))\bar u_\eps(x)u_\eps(x,t)\,dx .\label{msedp-eq-sysgen-asb3}
\end{align}

By Lemma  \ref{msedp-lem-estigen1} and \eqref{msedp-eq-sysgen-asb1}, we can check that  when $\eps \le \eps_1$ there exists positives  constants $c_1,C_1,\bar c_2,\bar C_2$  independent of $\eps$ such that for any positive smooth solutions $u_\epsilon$  to \eqref{msedp-eq-pert}--\eqref{msedp-eq-pert-bc} there exists $\bar t(u_\eps)$ so that   
\begin{equation}\label{msedp-eq-asb-esti1}
\hat c:=\frac{c_1}{\bar C_2}\le \lambda_\eps(t)\le \frac{C_1}{\bar c_2}=:\hat C \quad \text{for all }\quad t>\bar t.
\end{equation}
From the decomposition, by using  \eqref{msedp-eq-asb-esti1} and  Lemma \ref{msedp-lem-estigen1} we can also check that $h_\eps$ is smooth (i.e $C^{2,\alpha}(\O)$) and therefore belongs to $L^2(\O)$ for all times.

By plugging the decomposition of $u_\eps$ in \eqref{msedp-eq-sysgen-asb3} and using the definition of $\Psi_\eps$, we can check that 
\begin{equation}
\lambda^{\prime}_\eps(t)=\frac{\tilde \Psi_{\eps}(\bar u_\eps)\lambda_\eps(t)}{\nlp{\bar u_\eps}{2}{\O}^2}(1-\lambda_\eps^p(t)) + \r_1(t) +\r_2(t)
  \end{equation}
  where $\r_i$ are the following quantity:
  \begin{align}
  &\r_1(t):=\frac{1}{\nlp{\bar u_\eps}{2}{\O}^2}\int_{\O}[\Psi_{\eps}(x,\bar u_\eps)-\Psi_\eps(x,u_\eps)]\bar u_\eps(x)h_\eps(t,x)\,dx \\
  &\r_2(t):=\frac{\lambda_\eps(t)}{\nlp{\bar u_\eps}{2}{\O}^2}\int_{\O}\left(\sum_{k=1}^p\binom{k}{p}\lambda_\eps^{p-k}(t)\int_{\O}k_\eps(x,y)\bar u_\eps^{p-k}(y)h_\eps^{k}(t,y)\,dy\right)\bar u_\eps(x)^2\,dx 
  \end{align}

Next, we show  that 
\begin{lemma} \label{msedp-cla-energy-gen}
Let $p=1$ or $p=2$ then  there exists $\eps^*\le \min\{\eps_0,\eps_1\},$ so that for all $\eps\le\eps^*$ then any positive smooth solution $u_\eps$ of \eqref{msedp-eq-pert}--\eqref{msedp-eq-pert-bc} satisfies $$\lim_{t\to \infty}\oph{2}{\bar u_\eps}{h_{\eps}}(t))=0.$$  
 \end{lemma}

Assume the lemma holds true, then we can conclude  the proof of Theorem \ref{msedp-thm3} by arguing as follows. 
By combining Lemma \ref{msedp-lem-estigen1},  Lemma \ref{msedp-cla-energy-gen} and by using  H\"older's inequality, since $p=1$ or $2$ we see that  $\r_i(t) \to 0$ as $t \to +\infty$. Thus $\lambda_\eps(t)$ satisfies 
\begin{equation}
\lambda^{\prime}_\eps(t)=\frac{\tilde \Psi_{\eps}(\bar u_\eps)\lambda_\eps(t)}{\nlp{\bar u_\eps}{2}{\O}^2}(1+o(1)-\lambda_\eps^p(t)), 
  \end{equation}
The above ODE is of logistic type with a perturbation $o(1)\to 0$  with  a non negative initial datum.    Therefore, when $\eps\le \eps^*$    $\lambda_\eps(t)$ converges to  $1$ and we  conclude that  when $\eps\le \eps^*$ then any positive solution  $u_\eps$ to \eqref{msedp-eq-pert}--\eqref{msedp-eq-pert-bc} converges  to $\bar u_\eps$ almost everywhere.

\fdem

Let us now turn our attention to the proof of the Lemma \ref{msedp-cla-energy-gen}.

\dem{Proof of Lemma \ref{msedp-cla-energy-gen}:}

First, let us denote $\Gamma(t,x):=\Psi_\eps(x,\bar u_\eps)-\Psi_\eps(x,u_\eps)$.  By \eqref{msedp-eq-sysgen-asb2}  \eqref{msedp-eq-sysgen-asb3}  and by using Theorem \ref{msedp-thm-gen-id} we  achieve 

\begin{align*}
\frac{d\oph{2}{\bar u_\eps}{h_\eps}(t)}{dt}
&=-2 \int_{\O}\bar u_\eps^2\left(\nabla\left( \frac{h_\eps(t,x)} {\bar u_\eps(x)}\right)\right)^tA(x)\nabla\left( \frac{h_\eps} {\bar u_\eps}\right) +2\int_{\O}\Gamma(t,x) h_\eps(x) u_\eps(x)\, dx.
\end{align*}
Therefore using the definition of $\Psi_\eps$ and that $\bar u_\eps \perp h_\eps$ we have 
\begin{multline*}
\frac{d\oph{2}{\bar u_\eps}{h_\eps}(t)}{dt}= -2 \int_{\O}\bar u_\eps^2\left(\nabla\left( \frac{h_\eps(t,x)} {\bar u_\eps(x)}\right)\right)^tA(x)\nabla\left( \frac{h_\eps} {\bar u_\eps}\right)+2(\Psi_0(\bar u_\eps)-\Psi_0(u_\eps))\oph{2}{\bar u_\eps}{h_\eps}(t)\\ +2\eps\int_{\O}(\Psi_1(x,\bar u_\eps)-\Psi_1(x,u_\eps))  h_\eps^2(x)\, dx+2\eps\lambda_\eps\int_{\O}(\Psi_1(x,\bar u_\eps)-\Psi_1(x,u_\eps))  h_\eps(x)\bar u_\eps\, dx.
\end{multline*}
Let $\eps\le \min\{\eps_1,\eps_2\}$, by Lemma \ref{msedp-lem-estigen1} any stationary solution  $\bar u_\eps$ to \eqref{msedp-eq-pert}--\eqref{msedp-eq-pert-bc} is bounded in $L^p(\O)$ and for any positive solution $u_\eps$ to \eqref{msedp-eq-pert}--\eqref{msedp-eq-pert-bc} there exists $\bar t(u_\eps)$  so that for all times $t\ge \bar t$, $$c_p\le\nlp{u_\eps}{p}{\O}<C_p.$$ 
So for all times $t\ge \bar t$ we have
$$|\Psi_1(x,\bar u_\eps)-\Psi_1(x,u_\eps)|\le 2\|k_1\|_{\infty}\sup\{C_p,\bar C_p\}=: \kappa_1,$$
which implies that for $t\ge \bar t$
\begin{multline}
\frac{d\oph{2}{\bar u_\eps}{h_\eps}(t)}{dt}\le-2\int_{\O}\bar u_\eps^2(x)\left(\nabla\left( \frac{h_\eps(t,x)} {\bar u_\eps(x)}\right)\right)^tA(x)\nabla\left( \frac{h_\eps} {\bar u_\eps}\right)\,dx +2(\Psi_0(\bar u_\eps)-\Psi_0(u_\eps))\oph{2}{\bar u_\eps}{h_\eps}(t)\\+2\eps\kappa_1\oph{2}{\bar u_\eps}{h_\eps}(t) +2\eps\lambda_\eps\int_{\O}(\Psi_1(x,\bar u_\eps)-\Psi_1(x,u_\eps))  h_\eps(x)\bar u_\eps(x)\, dx.\label{msedp-eq-sysgen-asb4}
\end{multline}

By \eqref{msedp-eq-sysgen-asb1} \eqref{msedp-eq-sysgen-asb3}, using the definition of $\Psi_\eps$ we also have 
\begin{align*}
\frac{d}{dt}\oph{1}{\bar u_\eps}{u_\eps}(t)&=(\Psi_0(\bar u_\eps) -\Psi_0(u_\eps))\oph{1}{\bar u_\eps}{u_\eps}(t) +\eps\int_{\O}(\Psi_1(x,\bar u_\eps)-\Psi_1(x,u_\eps))  u_\eps(x)\bar u_\eps(x)\, dx.\\
&\ge (\Psi_0(\bar u_\eps)-\eps\kappa_1 -\Psi_0(u_\eps))\oph{1}{\bar u_\eps}{u_\eps}(t). 
\end{align*}
Since $\oph{1}{\bar u_\eps}{u_\eps}>0$ for all $t>0$,  we have
$$\frac{d \log(\oph{1}{\bar u_\eps}{u_\eps})}{dt}(t) \ge  (\Psi_0(\bar u_\eps)-\eps\kappa_1 -\Psi_0(u_\eps)),$$ 
which combined with \eqref{msedp-eq-sysgen-asb4} implies that for $t\ge \bar t$ 
  \begin{multline*}
\frac{d\oph{2}{\bar u_\eps}{h_\eps}(t)}{dt}\le-2\int_{\O}\bar u_\eps^2(x)\left(\nabla\left( \frac{h_\eps(t,x)} {\bar u_\eps(x)}\right)\right)^tA(x)\nabla\left( \frac{h_\eps} {\bar u_\eps}\right)dx+\left(\frac{d \log{\left(\oph{1}{\bar u_\eps}{u_\eps}\right)^2}}{dt}(t)\right) \oph{2}{\bar u_\eps}{h_\eps}(t)\\+4\eps\kappa_1\oph{2}{\bar u_\eps}{h_\eps}(t) +2\eps\lambda_\eps(t) \int_{\O}\Gamma_1(t,x)h_\eps(x)\bar u_\eps(x)\, dx.
  \end{multline*}
  where $\Gamma_1(t,x):=\Psi_1(x,\bar u_\eps)-\Psi_1(x,u_\eps)$.

Since $\eps \le \eps_1$,  by Lemma \ref{msedp-lem-esti-fci}, and by rearranging the terms in the above inequality we get  for $t\ge \bar t$ 
 \begin{multline}
\frac{d\oph{2}{\bar u_\eps}{h_\eps}(t)}{dt}-\oph{2}{\bar u_\eps}{h_\eps}(t)\frac{d\log{(\oph{1}{\bar u_\eps}{u_\eps})^2}}{dt}  (t)  \le (-\bar \rho+4\eps\kappa_1)\oph{2}{\bar u_\eps}{h_\eps}(t) \\+2\eps\lambda_\eps(t) \int_{\O}\Gamma_1(t,x)h_\eps(x)\bar u_\eps(x)\, dx.\label{msedp-eq-sysgen-asb6}
  \end{multline}

Now, we estimate the last term of the above inequality. 

\subsection*{Case $p=1$}
In this situation, by using the definition of $\Gamma_1$ and the Cauchy-Schwartz inequality we have 
$$
|\Gamma_1(t,x)|\le |1-\lambda_\eps|\nlto{\bar u_\eps} \sup_{x\in\O}\sqrt{\int_{\O}k_1(x,y)^2\,dy} +\nlto{h_\eps} \sup_{x\in\O}\sqrt{\int_{\O}k_1(x,y)^2\,dy}.
$$
Since $\nlto{v}=\sqrt{\oph{2}{\bar u_\eps}{v}}$, by the Cauchy-Schwartz inequality we achieve for $t\ge bar t$
\begin{align*}
\int_{\O}\Gamma_1(t,x)h_\eps(x)\bar u_\eps(x)&\le \kappa \sqrt{\oph{2}{\bar u_\eps}{\bar u_\eps}(t)} \sqrt{\oph{2}{\bar u_\eps}{h_\eps}(t)}\left[ |1-\lambda(t)| \sqrt{\oph{2}{\bar u_\eps}{\bar u_\eps}} +\sqrt{\oph{2}{\bar u_\eps}{h_\eps}(t)}\right],\\
&\le  \kappa \bar C_2 \sqrt{\oph{2}{\bar u_\eps}{h_\eps}(t)}\left[ |1-\lambda(t)| \bar C_2+\sqrt{\oph{2}{\bar u_\eps}{h_\eps}(t)}\right],
\end{align*}
where $\kappa:= \sup_{x\in\O}\sqrt{\int_{\O}k_1(x,y)^2\,dy}$.

  \subsection*{Case $p=2$}
In this situation, as above  by using the definition of $\Gamma_1$ and the Cauchy-Schwartz inequality, we see that

$$
|\Gamma_1(t,x)|\le |1-\lambda^2_\eps|\|k_1\|_\infty |\nlto{\bar u_\eps}^2+2 \lambda_\eps \|k_1\|_\infty\nlto{\bar  u_\eps}\nlto{h_\eps} + \|k_1\|_\infty \nlto{h_\eps(t)}^2.
$$

So we get for $t\ge \bar t$
\begin{align*}
\int_{\O}\Gamma_1(x)h_\eps(x)\bar u_\eps(x)&\le \kappa  \sqrt{\oph{2}{\bar u_\eps}{\bar u_\eps}(t)} \sqrt{\oph{2}{\bar u_\eps}{h_\eps}(t)}\left[ |1-\lambda^2_\eps| |\nlto{\bar u_\eps}^2+2 \lambda_\eps \nlto{\bar  u_\eps}\nlto{h_\eps} +  \nlto{h_\eps(t)}^2\right],\\
&\le  \kappa \bar C_2 \sqrt{\oph{2}{\bar u_\eps}{h_\eps}(t)}\left[ |1-\lambda^2(t)| \bar C_2^2+(\tilde C_2+\hat C)\sqrt{\oph{2}{\bar u_\eps}{h_\eps}(t)}\right],
\end{align*}
where $\kappa= \|k_1\|_{\infty}$.

\medskip

 In both case, we can see that there exists $\kappa_2$ and $\kappa_3$ independent of $\eps$, $\bar u_\eps$ and $u_\eps$ so that we have for   $t\ge \bar t$. 
\begin{equation}\label{msedp-eq-sysgen-asb5} 
\int_{\O}\Gamma_1(x)  h_\eps(x)\bar u_\eps(x)\, dx \le \kappa_2\sqrt{\oph{2}{\bar u_\eps}{h_\eps}(t)}\left[ |1-\lambda^p(t)|  +\kappa_3\sqrt{\oph{2}{\bar u_\eps}{h_\eps}(t)}\right].
\end{equation}

  By combining \eqref{msedp-eq-sysgen-asb5} and  \eqref{msedp-eq-sysgen-asb6}, we achieve  for $t\ge \bar t$ 
  \begin{multline}
\frac{d\oph{2}{\bar u_\eps}{h_\eps}(t)}{dt}-\oph{2}{\bar u_\eps}{h_\eps}(t)\frac{d\log{(\oph{1}{\bar u_\eps}{u_\eps})^2}}{dt}(t)  \le\left(-\bar \rho+\eps \kappa_5 \right)\oph{2}{\bar u_\eps}{h_\eps}(t) \\+\eps \kappa_4|1-\lambda^p(t)|\sqrt{\oph{2}{\bar u_\eps}{h_\eps}(t)},  \label{msedp-eq-sysgen-asb8}
  \end{multline}
  where   $\kappa_4:=2\hat C\kappa_2$  and $\kappa_5:= 2\hat C\kappa_2\kappa_3+4\kappa_1$ are positive constants independent of $\eps, u_\eps$ and $\bar u_\eps$. 

%
%

The proof now will follow several steps:
\subsubsection*{Step One:}
Since $\eps \le \eps_1$ by \eqref{msedp-eq-asb-esti1} we have $|1-\lambda^p(t)|\le \kappa_6$ for all $t>\bar t$, with $\kappa_6$ a universal constant independent of $\eps$.
We  claim that 
\begin{claim}
Let $\eps \le \eps_3:= \min( \eps_1,\eps_2:=\frac{\bar \rho}{2 \kappa_5}) $, then for all $u_\eps$ positive solution to \eqref{msedp-eq-pert}--\eqref{msedp-eq-pert-bc} there exists $\bar t'\ge \bar t$ so that for all $t\ge \bar t'$ we have 
$$\sqrt{\oph{2}{\bar u_\eps}{h_\eps}(t)}\le 2\eps \left(\frac{\hat C}{\hat c}\right)\frac{2\kappa_4\kappa_6}{\bar \rho}.$$
\end{claim}
\dem{Proof:}
Indeed for $\eps\le  \eps_3$ by  \eqref{msedp-eq-sysgen-asb8} for $t\ge \bar t$ we have 
\begin{equation}
\frac{d\oph{2}{\bar u_\eps}{h_\eps}(t)}{dt}-\oph{2}{\bar u_\eps}{h_\eps}(t)\frac{d}{dt}\log{(\oph{1}{\bar u_\eps}{u_\eps}(t))^2}  \le - \frac{\bar \rho}{2}\oph{2}{\bar u_\eps}{h_\eps}(t) +\eps \kappa_4\kappa_6\sqrt{\oph{2}{\bar u_\eps}{h_\eps}(t)}.   \label{msedp-eq-sysgen-asb9}
  \end{equation}

 From the above differential inequality we can check that there exists $t_0'>\bar t$ so that $$ \sqrt{\oph{2}{\bar u_\eps}{h_\eps}(t_0)}\le \frac{\eps 4\kappa_4\kappa_6}{\bar \rho}.$$  If not, then  $ \sqrt{\oph{2}{\bar u_\eps}{h_\eps}(t)}>\frac{\eps 4\kappa_4\kappa_6}{\bar \rho}$ for all $t>\bar t$  and by dividing  \eqref{msedp-eq-sysgen-asb9} by $\sqrt{\oph{2}{\bar u_\eps}{h_\eps}(t)}$  and 
 by rearranging the terms, we get  the inequality 
\begin{equation}
\sqrt{\oph{2}{\bar u_\eps}{h_\eps}(t)}\frac{d}{dt}\log{\left(\frac{\oph{2}{\bar u_\eps}{h_\eps}(t)}{\oph{1}{\bar u_\eps}{u_\eps}(t)^2}\right)}  \le -\frac{\bar \rho}{2}\sqrt{\oph{2}{\bar u_\eps}{h_\eps}(t)} +\eps \kappa_4\kappa_6< -\eps \kappa_4\kappa_6\qquad \forall t\ge \bar t.  \label{msedp-eq-sysgen-asb10}
  \end{equation}
  Thus $ F(t):=\log{\left(\frac{\oph{2}{\bar u_\eps}{h_\eps}(t)}{\oph{1}{\bar u_\eps}{u_\eps}(t)^2}\right)}$ is a decreasing function which  is bounded from below since $\lambda_\eps \le \hat C$.  Moreover $ \sqrt{\oph{2}{\bar u_\eps}{h_\eps}(t)}>\frac{\eps 4\kappa_4\kappa_6}{\bar \rho}$ for all $t>\bar t$. Therefore $F$ converges as $t$ tends to $+\infty$ and $\frac{dF}{dt} \to 0$. 
  Thus for $t$ large enough, we get the contradiction 
  $$-\frac{\eps \kappa_4\kappa_6}{2\bar \rho}\le \sqrt{\oph{2}{\bar u_\eps}{h_\eps}(t)}\frac{d}{dt}\log{\left(\frac{\oph{2}{\bar u_\eps}{h_\eps}(t)}{\oph{1}{\bar u_\eps}{u_\eps}(t)^2}\right)}  \le - \frac{\eps \kappa_4\kappa_6}{\bar \rho}.$$

Let $\Sigma$ be  the set  $\Sigma:=\left\{t>t_0' | \sqrt{\oph{2}{\bar u_\eps}{h_\eps}(t)}>  \frac{\eps 4\kappa_4\kappa_6}{\bar \rho}\right\}$.  
Assume that $\Sigma$ is non empty otherwise the claim is proved  since $\frac{\hat C}{\hat c}>1$.  Let us  denote   $t^*:=\inf\Sigma $. By construction, since $h_\eps$ is continuous we have $\sqrt{\oph{2}{\bar u_\eps}{h_\eps}(t^*)} = \frac{\eps 4\kappa_4\kappa_6}{\bar \rho}$. 

Again,  by dividing\eqref{msedp-eq-sysgen-asb9}  by $\sqrt{\oph{2}{\bar u_\eps}{h_\eps}(t)}$  and rearranging the terms,   we get for all   $t\in\Sigma$
\begin{equation}
\sqrt{\oph{2}{\bar u_\eps}{h_\eps}(t)}\frac{d}{dt}\log{\left(\frac{\oph{2}{\bar u_\eps}{h_\eps}(t)}{\oph{1}{\bar u_\eps}{u_\eps}(t)^2}\right)}  \le -\frac{\bar \rho}{4}\sqrt{\oph{2}{\bar u_\eps}{h_\eps}(t)} +\eps \kappa_4\kappa_6\le 0.  \label{msedp-eq-sysgen-asb11}
  \end{equation}

Thus $\log{\left(\frac{\oph{2}{\bar u_\eps}{h_\eps}(t)}{\oph{1}{\bar u_\eps}{u_\eps}(t)^2}\right)}$ is a decreasing function of $t$ for all $t \in \Sigma$. By arguing on each connected component of $\Sigma$ and by using Lemma \ref{msedp-lem-estigen1} we can check that  
for $t\ge t^*$ we have
$$\sqrt{\oph{2}{\bar u_\eps}{h_\eps}(t)}\le \frac{\hat C}{\hat c}\frac{\eps 4\kappa_4\kappa_6}{\bar \rho}.$$

Hence, since $\frac{\hat C}{\hat c}>1$ we get for all $t\ge t_0$, 
$$\sqrt{\oph{2}{\bar u_\eps}{h_\eps}(t)}\le \frac{ \hat C}{\hat c}\frac{\eps 4\kappa_4\kappa_6}{\bar \rho}.$$
\fdem

\subsubsection*{Step Two:}

Recall that $\lambda_\eps(t)$ satisfies 
\begin{equation}
\lambda_\eps^{\prime}(t)=\frac{\tilde \Psi_{\eps}(\bar u_\eps)\lambda_\eps(t)}{\nlp{\bar u_\eps}{2}{\O}^2}(1-\lambda_\eps^p(t)) + \r_1(t) +\r_2(t) \label{msedp-eq-sysgen-asb-edo}
  \end{equation}
  where $\r_i$ are the following quantity:
  \begin{align}
  &\r_1(t):=\frac{1}{\nlp{\bar u_\eps}{2}{\O}^2}\int_{\O}[\Psi_{\eps}(x,\bar u_\eps)-\Psi_\eps(x,u_\eps)]\bar u_\eps(x)h(t,x)\,dx \\
  &\r_2(t):=\frac{\lambda_\eps(t)}{\nlp{\bar u_\eps}{2}{\O}^2}\int_{\O}\left(\sum_{k=1}^p\binom{k}{p}\lambda_\eps^{p-k}(t)\int_{\O}k_\eps(x,y)\bar u_\eps^{p-k}(y)h^{k}(t,y)\,dy\right)\bar u_\eps(x)^2\,dx 
  \end{align}
Since $p=1$ or $p=2$ then by Lemma \ref{msedp-lem-estigen1} and  H\"older's inequality, we can see that  there exists $\kappa_7$  independent of $\eps, \bar u_\eps, u_\eps$ so that for all $t\ge \bar t$
\begin{equation}
|\r_1(t)+\r_2(t)|\le \kappa_7\frac{\tilde\Psi(\bar u_\eps)\lambda_\eps(t)}{\nlp{\bar u_\eps}{2}{\O}^2} \sqrt{\oph{2}{\bar u_\eps}{h_\eps}(t)}. \label{msedp-eq-sysgen-asb-R}
\end{equation}


Next, we define some constant quantities: 
\begin{align}
&\delta_0:=\frac{\hat C}{\hat c}\frac{4\kappa_4\kappa_6}{\bar \rho},\label{msedp-eq-sysgen-asb12}\\
&\eps^*:=\min\left\{\eps_3,\frac{\bar \rho \hat c}{16\kappa_4\kappa_7\hat C}, \frac{1}{4\kappa_7\delta_0}\right\},
\end{align}

By the previous step, we see that for $\eps \le \eps^*$  we have for any positive solution $u_\eps$ to \eqref{msedp-eq-pert}--\eqref{msedp-eq-pert-bc} there exists $\bar t'$ so that  for all $t\ge \bar t'$
$$ \sqrt{\oph{2}{\bar u_\eps}{h_\eps}(t)}\le \eps\delta_0.$$
We claim that 
\begin{claim}\label{msedp-cla-steptwo}
For  $\eps\le \eps^*$,   there exists $t_ {\eps\delta_0}\ge \bar t'$ such that for all $t\ge t_{\eps\delta_0}$
$$\sqrt{\oph{2}{\bar u_\eps}{h_\eps}(t)}\le \frac{\eps\delta_0}{2}.$$
 \end{claim}
 \dem{Proof:}	
First, we  can check that for $\eps\le\eps^*$  there exists $t^*$  so that for all $t\ge t^*$
$$|1-\lambda_\eps^p(t)|\le 2\eps\delta_0\kappa_7.$$

 Let $\lambda_{_{\pm \eps\delta_0\kappa_7}}\in C^1((\bar t',\infty),\R^+)$ be the solution of the ODE
 \begin{equation}
\lambda_{_{\pm \eps\delta_0\kappa_7}}'(t)=\frac{\tilde \Psi_{\eps}(\bar u_\eps)\lambda_{_{\pm \eps\delta_0\kappa_7}}}{\nlp{\bar u_\eps}{2}{\O}^2}(1\pm\eps\delta_0\kappa_7-\lambda_{_{\pm \eps\delta_0\kappa_7}}^p(t)),\quad \lambda_{_{\pm \eps\delta_0\kappa_7}}(\bar t')=\lambda_\eps(\bar t').\label{msedp-eq-sysgen-asb-edo1}
\end{equation}  
Since the above equation is of logistic type and $\lambda_{_{\pm \eps\delta_0\kappa_7}}(\bar t')>0$, $\lambda_{_{\pm \eps\delta_0\kappa_7}}(t)\to \bar \lambda_{\pm}$ as $t \to \infty$ where $\bar \lambda_{\pm}$ is the solution of the algebraic equation
$ 1\pm\eps\delta_0\kappa_7-\bar \lambda_{\pm}^p=0$.

By \eqref{msedp-eq-sysgen-asb-edo} and \eqref{msedp-eq-sysgen-asb-R},  we can check that $\lambda_\eps$ satisfies for $t \ge \bar t'$
\begin{align}
&\lambda_\eps^{\prime}(t)\ge\frac{\tilde \Psi_{\eps}(\bar u_\eps)\lambda_\eps(t)}{\nlp{\bar u_\eps}{2}{\O}^2}(1-\eps\delta_0\kappa_7 -\lambda_\eps^p(t)), \label{msedp-eq-sysgen-asb-edo2}\\
&\lambda_\eps^{\prime}(t)\le\frac{\tilde \Psi_{\eps}(\bar u_\eps)\lambda_\eps(t)}{\nlp{\bar u_\eps}{2}{\O}^2}(1+\eps\delta_0\kappa_7-\lambda_\eps^p(t)). \label{msedp-eq-sysgen-asb-edo3}
  \end{align}
 By the comparison principle, from \eqref{msedp-eq-sysgen-asb-edo1} \eqref{msedp-eq-sysgen-asb-edo2} and \eqref{msedp-eq-sysgen-asb-edo3}   we get $ \lambda_{-\eps\delta_0\kappa_7}(t)\le \lambda_\eps(t) \le \lambda_{+\eps\delta_0\kappa_7}(t)$ for all $t\ge \bar t'$. Thanks to the convergence of $\lambda_{\pm\eps\delta_0\kappa_7}(t)$ to $\bar \lambda_{\pm\eps\delta_0\kappa_7}$ and the monotone behaviour of  $\bar \lambda_{\pm\eps\delta_0\kappa_7}$ with respect to $\eps$ we get  
$$\bar \lambda_{-2\eps\delta_0\kappa_7}\le \lambda_\eps(t) \le \bar \lambda_{+2\eps\delta_0\kappa_7} \quad\text{for}\quad t\ge t^*,$$
for some $t^*\ge \bar t'$.
Therefore, for $t\ge t^*$ we have 
$$|1-\lambda_\eps^p(t)|\le 2\eps\delta_0 \kappa_7.$$

From the latter estimate,  since $\eps \le \eps_3$  we deduce from  \eqref{msedp-eq-sysgen-asb8} that for $t\ge t^*$

\begin{equation*}
\frac{d\oph{2}{\bar u_\eps}{h_\eps}(t)}{dt}-\oph{2}{\bar u_\eps}{h_\eps}(t)\frac{d}{dt}\log{(\oph{1}{\bar u_\eps}{u_\eps}(t))^2}  \le -\frac{\bar \rho}{2}\oph{2}{\bar u_\eps}{h_\eps}(t) +2\eps^2 \kappa_4\kappa_7 \delta_0\sqrt{\oph{2}{\bar u_\eps}{h_\eps}(t)}.  
  \end{equation*}
  By following the argumentation of Step one,  we can show that there exists $ t_{\eps\delta_0}\ge t^*$ such that for   $t\ge t_{\eps\delta_0}$ we have 
$$\sqrt{\oph{2}{\bar u_\eps}{h_\eps}(t)}\le \frac{8\eps \hat C \kappa_4 \kappa_7}{\hat c \bar \rho} \eps\delta_0,$$ 
which thanks to $\eps\le  \frac{\bar \rho \hat c}{16\kappa_4\kappa_7\hat C}$ leads to 
  
$$\sqrt{\oph{2}{\bar u_\eps}{h_\eps}(t)}\le \frac{\eps\delta_0}{2}.$$  

\fdem

\subsubsection*{Step Three:}
Since for all $t\ge t_{\eps\delta_0}$, 
$$\sqrt{\oph{2}{\bar u_\eps}{h_\eps}(t)}\le \frac{\eps\delta_0}{2},$$
 by arguing as in the proof of Claim \ref{msedp-cla-steptwo}, we see that there exists $t_{\eps\frac{\delta_0}{2}}$ so that for all $t\ge t_{\eps\frac{\delta_0}{2}}$
$$\sqrt{\oph{2}{\bar u_\eps}{h_\eps}(t)}\le \frac{\eps\delta_0}{4}.$$ 

By reproducing inductively the above argumentation,  we can construct a sequence  $(t_{n})_{n\in \N}$ so that 
for all $t\ge t_n$ we have 
$$ \sqrt{\oph{2}{\bar u_\eps}{h_\eps}(t)}\le \frac{\eps\delta_0}{2^n}.$$
Hence, when $\eps \le \eps^*$ we deduce that  
$$\lim_{t\to \infty}\oph{2}{\bar u_{\eps}}{h_\eps}(t) \to 0.$$
   
\fdem

\bigskip

\noindent \textbf{Acknowledgements.}  The author thanks the members of the INRIA project: ERBACE,  for early discussion on this subject. 
The author wants also to thanks Professor Raoul for interesting discussions on these topic.


\appendix

\appendix
\section{Existence of a positive solution}

In this appendix, we present a  construction of a smooth positive solution of \eqref{msedp-eq-intro} 
 The construction  is rather simple and follows some of the ideas used in \cite{Calsina2012}. 
 First, let $p\ge 1$ be fixed and let us regularised $u_0$ by a smooth mollifier $\rho_{\eps}$ and consider the solution  of \eqref{msedp-eq-intro}--\eqref{msedp-eq-intro-ci} with initial condition $u_{\eps,_0}:=\rho_\eps\star u_0$ instead of $u_0$.
 
Now we introduce the following sequence of function $(u_n(x,t))_{n\in \N}$ where $u_{n}$ is defined recursively by
$u_0(x,t)=u_0(x)$ and for $n\ge 0$, $u_{n+1}$ is the solution of
   \begin{align}
&\frac{\partial u_{n+1}}{\partial t}=\nabla\cdot(A(x)\nabla u_{n+1} )+ u_{n+1}\left(r(x)-\int_{\O}K(x,y)|u_n|^p(t,y)\,dy\right)\quad \text{ in } \quad \R^+\times\O \label{msedp-eq-approx}\\
&\frac{\partial u_{n+1}}{\partial n}(t,x)=0 \quad \text{ in } \quad \R^+\times\partial\O\label{msedp-eq-approx-bc}\\
&u_{n+1}(x,0)=u_{\eps,0}(x)\quad \text{ in } \quad \O.\label{msedp-eq-approx-ci}
\end{align} 
Since by assumption $u_{\eps,0} \in C^{\infty}(\O)$, $(u_n)_{n\in\N}$ is well defined from the standard parabolic theory see \cite{Brezis2010,Evans1998}.
Moreover since $u_{\eps,0}\ge 0$  and  $0$ is a sub-solution of  the problem \eqref{msedp-eq-approx}-- \eqref{msedp-eq-approx-ci} for each $n$, by the parabolic strong maximum principle  we deduce  that $u_n(x,t)>0$ for all $n, x$ and $t>0$.

 Now since  $u_n$and $K$ are  non-negative functions, for all $n\ge 0, $ $u_{n+1}$ is a subsolution of the linear problem:
  \begin{align}
&\frac{\partial v}{\partial t}=\nabla\cdot(A(x)\nabla v )+ r(x)v\quad \text{ in } \quad \R^+\times\O \label{msedp-eq-cp}\\
&\frac{\partial v}{\partial n}(t,x)=0 \quad \text{ in } \quad \R^+\times\partial\O \label{msedp-eq-cp-bc}\\
&v(x,0)=u_0^\eps(x)\quad \text{ in } \quad \O.\label{msedp-eq-cp-ci}
\end{align} 
 and by the parabolic maximum principle,  we have $u_n\le v\le \|u_0^{\eps}\|_{\infty}e^{\|r\|_{\infty} t}$ in $\R^+\times \O$ for all $n$. Therefore from the standard Schauder parabolic \textit{a priori } estimates, we deduce that $(u_n)_{n\in \N}$ is uniformly bounded in $C^{1,\alpha}((0,T),C^{2,\beta}(\O))$ for each $T>0$. Thus  by diagonal extraction, there exists a subsequence $(u_{n_k})_{k\in\N}$ which converges to a solution $u(x,t)\ge 0$ of \eqref{msedp-eq-intro}--\eqref{msedp-eq-intro-ci} with initial condition $u_{\eps,0}$.


Let us now take the limit $\eps \to 0$. 
By multiplying \eqref{msedp-eq-intro}  by $\phi_1$ and integrate it over $\O$ we have

$$\frac{d}{dt}\left(\int_{\O}u_\eps(t,x)\phi_1(x)\,dx\right)= -\lambda_1 \int_{\O}u_\eps \phi_1- \int_{\O\times\O}K(x,y)\phi_1(x)u_\eps(t,x)u^p_\eps(t,y)\,dydx.$$
Since  $u^\eps, \phi_1$ and $K(x,y)$ are positives in $\bar \O$ it follows that 
$$\frac{d}{dt}\left(\int_{\O}u_\eps(t,x)\phi_1(x)\,dx\right)\le  \int_{\O}u_\eps \phi_1\left(-\lambda_1-  C_0\int_{\O}u_\eps \phi_1\right),$$
for some positive constant $C_0$ which depends only on $\phi_1$ and $K$. Thanks to the logistic character of the above inequality,  we deduce that $\nlp{u_\eps}{1}{\O}$ is bounded uniformly in time independently of  $\eps$. 
By using Theorem \ref{msedp-thm-gen-id} and Remark \ref{msedp-rem-id} with $H(s):s\mapsto s^p$ and $\phi_1$, it follows that 
 \begin{multline*}
 \frac{d \oph{p}{\phi_1}{u_\eps}(t)}{dt}\le - p(p-1)\int_{\O}\left(\frac{u_\eps(t,x)}{\phi_1(x)} \right)^{p-2}\phi_1^2 \left(\nabla\left( \frac{u_\eps(t,x)} {\phi_1(x)}\right)\right)^tA(x)\nabla\left( \frac{u_\eps(t,x)} {\phi_1(x)}\right)\,dx \\ +  p \int_{\O} \phi_1^2(x) \left(\frac{u_\eps}{\phi_1}(t,x)\right)^{p}\left[-\lambda_1-\int_{\O}K(x,y)u^p_\eps(t,y)\,dy\right]\,dx.
 \end{multline*}
As above since  $u^\eps, \phi_1$ and $K(x,y)$ are positives in $\bar \O$ it follows that 
$$ \frac{d \oph{p}{\phi_1}{u_\eps}(t)}{dt}\le   C_1 \oph{p}{\phi_1}{u_\eps}(t)\left[-\lambda_1-C_2\oph{p}{\phi_1}{u_\eps}(t)\right],$$
for some positive constants $C_1$ and $C_2$ which depends only on $\phi_1$ and $K$.
Thus $\nlp{u_\eps}{p}{\O}$ is bounded uniformly with respect to $\eps$.
Since the coefficient of the parabolic PDE are bounded in $L^\infty$ independently of $\eps$, by standard parabolic $L^p$ estimates \cite{Wu2006}, it follows that for all $T>0$, $u_\eps$ is  bounded independently of $\eps$ in   $W^{1,2,1}((0,T)\times \O)\cap W^{1,1,1}_{0}(0,T)\times \O)$, where for $p\ge 1$ $W^{1,2,p}$ and $ W^{1,1,p}_{0}$ denote the Sobolev space 
  \begin{align*}
  &W^{1,2,p}:=\{u\in L^p((0,T)\times\O)\,|\, \partial_t u, \nabla u, \partial_{ij}u \in L^p((0,T)\times\O) \},\\
  &W^{1,1,p}_0:=\{u\in L^p((0,T)\times\O), \partial_nu=0 \text{ on } \partial \O\,|\, \partial_t u, \nabla u \in L^p((0,T)\times\O) \}.
\end{align*}   
 
By a standard bootstrap argument using the Parabolic regularity, we see that for each $T>0$, $(u_\eps)$ is  bounded in $C^{1,\alpha}((0,T),C^{2,\beta}(\O))$ independently of $\eps$. Thus by diagonal extraction, there exists a subsequence $(u_{\eps_{n_k}})_{k\in\N}$ which converges to a smooth solution $u(x,t)\ge 0$ of \eqref{msedp-eq-intro}--\eqref{msedp-eq-intro-ci} with initial condition $u_0$.    

%
 
\section*{}
\bibliographystyle{amsplain}
\bibliography{msedp.bib}
\end{document}